\documentclass[a4paper]{amsart}
\usepackage[margin=3cm]{geometry}
\usepackage{amsthm}
\usepackage{enumerate}

\theoremstyle{plain}
\newtheorem{thm}{Theorem}[section]
\newtheorem*{thm*}{Theorem}
\newtheorem{prop}{Proposition}[section]
\newtheorem*{prop*}{Proposition}
\newtheorem{lemma}{Lemma}[section]
\newtheorem*{lemma*}{Lemma}

\newtheorem*{coro*}{Corollary}

\theoremstyle{definition}
\newtheorem{dfn}{Definition}[section]
\newtheorem*{dfn*}{Definition}
\newtheorem{rem}{Remark}[section]
\newtheorem*{rem*}{Remark} 
\newtheorem{ex}{Example}[section]
\newtheorem*{ex*}{Example}

\theoremstyle{remark}
\newtheorem*{rems*}{Remarks}
\newtheorem*{fact*}{Fact}

\newtheorem*{exc*}{Exercise}

\newcommand{\A}{\mathbb A}
\newcommand{\C}{\mathbb{C}}
\newcommand{\R}{\mathbb R}
\newcommand{\Z}{\mathbb Z}
\newcommand{\Q}{\mathbb Q}
\newcommand{\Gal}{\operatorname{Gal}}
\newcommand{\Hom}{\operatorname{Hom}}
\newcommand{\Aut}{\operatorname{Aut}}
\newcommand{\End}{\operatorname{End}}
\newcommand{\Res}{\operatorname{Res}}
\newcommand{\GL}{\operatorname{GL}}
\newcommand{\SL}{\operatorname{SL}}
\newcommand{\Ad}{\operatorname{Ad}}
\newcommand{\ad}{\operatorname{ad}}
\newcommand{\ab}{\operatorname{ab}}
\newcommand{\Gm}[1]{\mathbb{G}_{\operatorname{m},#1}}
\newcommand{\id}{\operatorname{id}}
\newcommand{\diag}{\operatorname{diag}}
\newcommand{\Sh}{\operatorname{Sh}}
\newcommand{\can}{\operatorname{can}}
\newcommand{\der}{\operatorname{der}}
\newcommand{\Lie}{\operatorname{Lie}}
\newcommand{\Spec}{\operatorname{Spec}}
\newcommand{\SO}{\operatorname{SO}}
\newcommand{\Inn}{\operatorname{Inn}}
\newcommand{\inte}{\operatorname{int}}

\numberwithin{equation}{section}

\begin{document} 

\title{Complex conjugation and Shimura varieties}
\author{Don Blasius, Lucio Guerberoff}

\maketitle

\begin{abstract} In this paper we study the action of complex conjugation on Shimura varieties and the problem of descending these to the maximal totally real field of the reflex field. We prove the existence of such descent for many Shimura varieties whose associated adjoint group has certain factors of type A or D. This includes a large family of Shimura varieties of abelian type. Our considerations and constructions are carried out purely at the level of Shimura data and group theory.
\end{abstract}

\section{Introduction}
The goal of this paper is to analyze some aspects of complex conjugation acting on Shimura varieties. This topic has been studied for a long time by several authors, notably Shimura, Deligne, Langlands, Milne, Shih, and more recently, Taylor. In general, given a Shimura variety $\Sh(G,X)$ defined by a Shimura datum $(G,X)$, and any automorphism $\alpha$ of $\C$, Langlands' conjectured (\cite{langlands}) that the conjugate variety $\alpha\Sh(G,X)=\Sh(G,X)\times_{\C,\alpha}\C$ can be realized as a Shimura variety $\Sh({}^{\alpha}G,{}^{\alpha}X)$ for a very explicit pair $({}^{\alpha}G,{}^{\alpha}X)$. This has been proved by Milne in \cite{milne83} (see also \cite{borovoiconcerning}, \cite{borovoipoints} and \cite{milnedescent}). The case of $\alpha=c$ (complex conjugation) has, among other properties, the particularity that the pair $({}^{c}G,{}^{c}X)$ is very concrete. Namely, it can be identified with $(G,\overline X)$, where $\overline X$ is obtained by composing the elements of $x$ with complex conjugation on the Deligne torus $\mathbb{S}$. This simple description is hard to find in the literature, and hence we include a proof of how it's deduced from the general constructions. 

Assuming a few standard extra conditions on the Shimura datum $(G,X)$, the reflex field $E$ can be seen to be either totally real or a CM field. The Shimura variety has a canonical model $\Sh(G,X)_{E}$ over $E$, and the Hecke operators are defined over $E$ as well. In this paper we investigate descent of these varieties to the maximal totally real subfield $E^{+}$ of $E$. The existence of such descent can be seen as a nice generalization of the useful fact that the field obtained by adjoining to $\Q$ the $j$-invariant of an order in an imaginary quadratic field has a real embedding. From now on, assume that $E$ is CM. We show in many cases that $\Sh(G,X)$ has a model over $E^{+}$. Although the Hecke operators are not defined over $E^{+}$, they can nevertheless be characterized. The general framework for constructing such models comes from the construction of descent data arising from automorphisms of $G$ of order $2$ taking $X$ to $\overline X$. Using the classification of (adjoint) Shimura data in terms of special nodes on Dynkin diagrams, our aim is to construct an involution of $G$ that induces the opposition involution on the based root datum (or the Dynkin diagram). The construction we make follows from the classification of semisimple groups. The groups $G$ which we will work with are, roughly speaking, those for which the simple factors of $G^{\ad}$ are of classical type $A$ or $D$, and satisfy an extra condition on the hermitian or skew-hermitian space defining them (see Definitions \ref{dfn:strhermitian} and \ref{dfn:strskewhermitian}). For example, a factor of type $A$ is attached to a hermitian space over a central division algebra $D$ over a CM field $K$ endowed with an involution of the second kind $J$. We show that if there exists an opposition involution on these groups, then $D$ must be either $K$ or a quaternion division algebra, and the involution $J$ is easily described. We carry out the construction of involutions if we assume the aforementioned extra condition, which in this case amounts to the existence of a basis of the underlying vector space such that the matrix of the hermitian form is diagonal with entries in $K$. In the quaternion algebra case, we can write $D=D_{0}\otimes_{F}K$, where $F$ is the maximal totally real subfield of $K$, and $D_{0}$ is a quaternion division algebra over $F$. We assume furthermore in this case that if $D_{0,v}$ is not split for an embedding $v:F\hookrightarrow\R$, then the corresponding factor of $G_{\R}^{\ad}$ is compact. If $D=K$, the conditions in Definition \ref{dfn:strhermitian} are automatically satisfied. For factors of type $D$, there is a similar scenario, although we only restrict to groups of type $D^{\mathbb{H}}$ as in the Appendix of \cite{milneshih}. This encompasses a large family of Shimura varieties of abelian type. We stress here that our methods are group-theoretic and we work purely at the level of Shimura data, without making use of any moduli interpretation. In particular, we could also include factors of type $E_{6}$ (which is the only other type apart from $A$ or $D^{\mathbb{H}}$ that contributes to the reflex field being CM instead of totally real) if they are concrete enough to construct involutions on them. 

Let us describe the organization of the paper and outline the main argument. In Section \ref{sec:shimura}, we start by recalling the general formalism of conjugation of Shimura varieties by an arbitrary automorphism of $\C$, we study the special case of complex conjugation explicitly, and prove in this case that the conjugate Shimura datum is $(G,\overline{X})$, where $\overline{X}$ is the complex conjugate conjugacy class of $X$. We show (Theorem \ref{thmdescent}) that if $(G,X)$ is a Shimura datum and $\theta:G\to G$ is an involution such that $\theta(X)=\overline{X}$, then $\theta$ induces an isomorphism of algebraic varieties from the complex conjugate $c\Sh(G,X)$ to $\Sh(G,X)$, defined over the reflex field $E$, that constitutes a descent datum from $E$ to $E^{+}$.

In Section \ref{sec:opposition}, we recall some basic facts about root data and opposition involutions, and in Proposition \ref{mainprop}, we lay the ground for the prototype of involutions $\theta:G\to G$ that we will construct. Roughly speaking, suppose that $T\subset G$ is a maximal torus of $G$, and $x\in X$ factors through $T_{\R}$. If $\theta:G\to G$ is an involution that preserves $T_{\R}$ and induces complex conjugation on the group of characters $X^{*}(T)$, then $\theta(x)=\overline{x}$ and thus $\theta(X)=\overline{X}$. This is basically the type of involutions that we will construct, with some slight changes. Since we will make use of the explicit classification of semisimple groups, we need to work with either $G^{\der}$ or $G^{\ad}$. We let $G_{i}$ be the almost simple factors of $G^{\der}$, and $\widetilde G_{i}$ be their simply connected covers, so that $\widetilde G_{i}=\Res_{F_{i}/\Q}H_{i}$, for certain groups $H_{i}$ which are absolutely almost simple, simply connected, over a totally real field $F_{i}$. We recall the classification of these groups in Section \ref{sec:semisimple}, where we also construct opposition involutions on them preserving specific maximal tori $S_{i}$ and inducing complex conjugation on their characters (for non-compact places $v$ of $F_{i}$). We only do this for groups of type $A$ or $D^{\mathbb{H}}$. These, together with type $E_{6}$, are the only ones that give a CM reflex field, as opposed to totally real. Furthermore, as noted above, we impose some extra conditions in order to construct the involutions. From the tori $S_{i}$, we get maximal tori $T'\subset G^{\der}$ and $T\subset G$, and an opposition involution $\theta':G^{\der}\to G^{\der}$ preserving $T'$. As shown in Proposition \ref{mainprop}, $\theta'$ extends uniquely to an involution on $G$. To show that $\theta(X)=\overline{X}$, we need to relate in some way the choice of our tori $S_{i}$, which is a priori unrelated to the Shimura datum, to the conjugacy class $X$. In Section \ref{sec:invoshimura}, we show that there always exists $x\in X$ such that $x^{\ad}$ factors through the image of $T_{\R}$ in $G_{\R}^{\ad}$. This is all we need for Proposition \ref{mainprop}. In Theorem \ref{thmforstrong}, we state the existence of descent datum for Shimura varieties defined by groups $(G,X)$ such that the simple factors of $G^{\ad}$ are of the type described in Section \ref{sec:semisimple}. We call these {\em strongly of type} $(AD^{\mathbb{H}})$. Finally, we also note that involutions inducing the desired descent datum on $\Sh(G,X)$ can be constructed whenever $G$ is adjoint and there exists an opposition involution $\theta:G\to G$. This is always the case if $G$ is adjoint and quasi-split, for example. 

The existence of the involutions constructed in this paper should have interesting applications, which will be explored in the future, for example, in the setting of integral models and the zeta function problem, and periods of automorphic forms. 

The first author thanks the Erwin Schr\"odinger International Institute for Mathematics and Physics for its support during a visit when a part of the work for this project was undertaken. A substantial part of this work was carried out while the second author was a guest at the Max Planck Institute for Mathematics in Bonn, Germany. It is a pleasure to thank the Institute for its hospitality and the excellent working conditions.

\subsection*{Notation and conventions}

We fix an algebraic closure $\C$ of the real numbers $\R$, a choice of $i=\sqrt{-1}$, and we let $\overline\Q$ denote the algebraic closure of $\Q$ in $\C$. We let $c\in\Gal(\C/\R)$ denote complex conjugation on $\C$, and we use the same letter to denote its restriction to $\overline\Q$. Sometimes we also write $c(z)=\bar z$ for $z\in\C$.

Let $k$ be a field. By a variety over $k$ we will mean a geometrically reduced scheme of finite type over $k$. We let $\Gm{k}$ denote the usual multiplicative group over $k$. For any algebraic group $G$ over $k$, we let $\Lie(G)$ denote its Lie algebra and $\Ad:G\to\GL_{\Lie(G)}$ the adjoint representation. For us, reductive group will always include connectedness in the definition. If $G$ is reductive, we let $G^{\ad}$ (resp. $G^{\der}$) denote its adjoint group $G/Z(G)$ (resp. its derived subgroup), where $Z(G)$ is the center of $G$. If $T\subset G$ is a torus, we denote by $T^{\ad}$ the image of $T$ under the projection $G\to G^{\ad}$. For any commutative group scheme $G$, we denote by $i_{G}:G\to G$ the map $g\mapsto g^{-1}$.

We denote by $\A$ (resp. $\A_f$) the ring of ad\`eles of $\Q$ (resp. finite ad\`eles). A CM field $K$ is a totally imaginary quadratic extension of a totally real field $F$. 

We let $\mathbb{S}=R_{\C/\R}\Gm{\C}$. We denote by $c=c_{\mathbb{S}}$ the algebraic automorphism of $\mathbb{S}$ induced by complex conjugation. For any $\R$-algebra $A$, this is $c\otimes_{\R}1_A:(\C\otimes_{\R}A)^\times\to(\C\otimes_{\R}A)^\times$ on the points of $\mathbb{S}(A)$.This is often denoted by $z\mapsto \bar z$, and on complex points it should not be confused with the other complex conjugation $1_{\mathbb{S}}\otimes c$ on $\mathbb{S}(\C)=(\C\otimes_{\R}\C)^{\times}$ on the second coordinate.

An involution of a group is an automorphism of order $2$, whereas an involution of a ring is an antiautomorphism of order $2$. This should not cause any confusion.

We will denote by $\mathbb{H}$ the non-split quaternion algebra over $\R$, identified with the set of matrices of the form 
\[ \left(\begin{array}{cc} x & y \\ -\overline y & \overline x\end{array}\right)\] 
in $M_{2}(\C)$.

\section{Shimura varieties, conjugation and descent}\label{sec:shimura}
We will first review some basic facts about Shimura varieties and conjugation by an automorphism of $\C$, specializing to the case of complex conjugation. Then we set up our descent problem, describe some general considerations about reflex fields and Dynkin diagrams, and explain how to construct descent data based on involutions of a Shimura datum.

\subsection{Shimura varieties}\label{subsec:notshimura} A Shimura datum $(G,X)$ will be understood in the sense of Deligne's axioms (2.1.1.1-3, \cite{delignecorvallis}). We will assume moreover that the connected component $Z^{0}$ of the center $Z$ of $G$ splits over a CM field. For a compact open subgroup $K\subset G(\A_f)$, we put $\Sh_K(G,X)(\C)=G(\Q)\backslash X\times G(\A_f)/K$. For $K$ sufficiently small (which we assume from now on), this complex analytic space is smooth and is equal to the complex points of a complex quasi-projective variety $\Sh_K(G,X)_{\C}$. Let $E=E(G,X)\subset\C$ be the reflex field of $(G,X)$; under our hypotheses, this is contained in a CM field, and thus it's either a CM field or a totally real field. In any case, we let $E^+$ be the maximal totally real subfield of $E$. The variety $\Sh_K(G,X)_{\C}$ admits a canonical model over $E$, denoted by $\Sh_K(G,X)_{E}$. We use the same notation for the pro-objects $\Sh(G,X)(\C)$, $\Sh(G,X)_{\C}$ and $\Sh(G,X)_E$. We denote by $w_X:\Gm{\R}\to G_{\R}$ the composition of $x\in X$ with the weight morphism $w:\Gm{\R}\to\mathbb{S}$, for some (or any) $x\in X$, and call it the {\em weight morphism} of $(G,X)$. For $x\in X$, we let $\mu_{x}:\Gm{\C}\to G_{\C}$ be the map given by $\mu_{x}(z)=x_{\C}(z,1)$, under the identification of $\mathbb{S}_{\C}\cong\Gm{\C}\times\Gm{\C}$ given by $(z\otimes a)\mapsto (za,\bar za)$.

We will fix the following notation once and for all. Let $G^{\der}$ be the derived group of $G$, $G^{\ab}=G/G^{\der}$ (a torus), $G^{\ad}$ be the adjoint group of $G$ and $p:G\to G^{\ad}$ be the projection onto $G^{\ad}$. The natural isogeny $Z^{0}\times G^{\der}\to G$ and the projection $G\to G^{\ab}$ define an isogeny $Z^{0}\to G^{\ab}$. Let $G_{1},\dots,G_{r}$ be the almost simple factors of $G^{\der}$ over $\Q$, and let $\widetilde G_{i}\to G_{i}$ be their simply connected covers. We can write $\widetilde G_{i}=\Res_{F_{i}/\Q}H_{i}$, where the fields $F_{i}$ are totally real and the groups $H_{i}$ are simply connected, absolutely almost simple over $F_{i}$. For each embedding $v\in I_{i}=\Hom(F_{i},\C)$, we have groups $H_{i,v}=H_{i}\otimes_{F_{i},v}\R$, and for a fixed $i=1,\dots,r$, all these groups have the same Dynkin type $D_{i}$, which will be called the Dynkin type of $\widetilde G_{i}$ (or of $G_{i}$ or $H_{i}$). We let $I_{i,c}=\{v\in I_{i}\mid H_{i,v}^{\ad}(\R)\text{ is compact}\}$ and we let $I_{i,nc}$ be its complement in $I_{i}$, which must be non-empty if $H_{i}$ is non-trivial.  We also have that $G^{\ad}$ is the direct product of the $G_{i}^{\ad}=\Res_{F_{i}/\Q}H_{i}^{\ad}$, and $G^{\ad}_{\R}$ is the direct product of the $H_{i,v}^{\ad}$ for $i=1,\dots,r$ and $v\in I_{i}$. Let $X^{\ad}$ be the $G^{\ad}(\R)$-conjugacy class containing $p_{\R}(X)$, and write $X^{\ad}=\prod_{i,v}X_{i,v}$ with $X_{i,v}$ an $H_{i,v}^{\ad}(\R)$-conjugacy class of morphisms $\mathbb{S}\to H_{i,v}^{\ad}$. For each $i$ and each $v\in I_{i,nc}$, there is a special node $s_{i,v}$ in the Dynkin diagram $D_{i,v}$ of $H_{i,v}$ attached to $X_{i,v}$, which uniquely determines $X_{i,v}$ as a conjugacy class with target $H_{i,v}^{\ad}$ (in the sense that if $Y$ is an $H_{i,v}^{\ad}(\R)$-conjugacy class satisfying Deligne's axioms, for which its associated special node is $s_{i,v}$, then $Y=X_{i,v}$; see \cite{delignecorvallis}, 1.2.6).

\subsection{Conjugation} For the general properties of conjugation of Shimura varieties, we mainly follow \cite{dmos}; see also \cite{milne} and \cite{langlands}. Let $(G,X)$ be a Shimura datum. A special pair $(T,x)$ consists of a maximal torus $T\subset G$ and a point $x\in X$ factoring through $T_{\R}$. Fix $x\in X$ a special point, and let $\sigma\in\Aut(\C)$. We denote by $({}^{\sigma,x}G,{}^{\sigma,x}X)$ the conjugate Shimura datum, so that there exists an isomorphism $\varphi_{\sigma,x}:\sigma\Sh(G,X)_{\C}=\Sh(G,X)_{\C}\times_{\C,\sigma}\C\simeq\Sh({}^{\sigma,x}G,{}^{\sigma,x}X)_{\C}$, unique with certain natural properties (Theorem II.4.2, \cite{milne}). Choosing a different special point gives canonically isomorphic results (Proposition II.4.3, \cite{milne}). The reflex field of $({}^{\sigma,x}G,{}^{\sigma,x}G)$ is $\sigma(E)$, and $\sigma\Sh(G,X)_{E}=\Sh(G,X)_{E}\times_{E,\sigma}\sigma(E)$ is the canonical model of $\Sh({}^{\sigma,x}G,{}^{\sigma,x}G)_{\C}$ over $\sigma(E)$. All of this also works at finite level: if $K\subset G(\A_f)$ is compact open, $\varphi_{\sigma,x}$ sends $\sigma\Sh_K(G,X)_{\C}$ to $\Sh_{{}^{\sigma,x}K}({}^{\sigma,x}G,{}^{\sigma,x}X)_{\C}$ (same thing replacing $\C$ by $E$ and $\sigma(E)$), where ${}^{\sigma,x}K\subset{}^{\sigma,x}G(\A_f)$ is explicit (see below).

We are interested mainly on the case $\sigma=c$, but nevertheless it will be useful to recall the general construction of $({}^{\sigma,x}G,{}^{\sigma,x}X)$. Let $\mathfrak{S}$ be the (connected) Serre group. This can be defined as the group of automorphisms of the forgetful fibre functor from the Tannakian category of CM $\Q$-Hodge structures to the category of finite dimensional $\Q$-vector spaces. (Here  a $\Q$-Hodge structure is a $\Q$-vector space $V$ such that $V\otimes \C$ is endowed with a Hodge structure; the structure is CM if the  algebra of elements of $\End(V)$ which induce morphisms of Hodge structure contains a commutative semisimple subalgebra of dimension $\dim _{\Q} (V)$).  Let $\mathfrak{T}$ denote the Taniyama group, defined here as the group of automorphisms of the Betti fibre functor in Deligne's Tannakian category of CM motives for absolute Hodge cycles over $\Q$; this is the Tannakian category generated by Artin motives and by the cohomology of abelian varieties over $\Q$ which are potentially CM. These are pro-algebraic groups, and there is a natural exact sequence
\[ 1\to\mathfrak{S}\to\mathfrak{T}\stackrel{\pi}{\to}\Gal({\overline\Q}/\Q)\to1,\]
where the second arrow corresponds to the functor taking a CM motive $M$ to its CM Hodge structure $H_B(M)$, and $\pi$ corresponds to the natural inclusion of the category of Artin motives into the category of CM motives. The group $\Gal({\overline\Q}/{\Q})$ is to be considered as the pro-algebraic group given by the inverse limit of the finite constant groups $\Gal(L/\Q)$, for $L\subset\C$ a finite Galois extension of $\Q$. There is a continuous section of $\pi$ over $\A_f$ denoted by $sp:\Gal({\overline\Q}/\Q)\to\mathfrak{T}(\A_f)$. For a motive $M$, $sp(\sigma)$ corresponds to the automorphism of $H_B(M)\otimes_{\Q}\A_f$ obtained  from the Galois action of $\sigma$ on \'etale cohomology using the comparison isomorphism. In the case of complex conjugation $c$,  $sp(c)=F\otimes_{\Q}1_{\A_f}$ for a unique  $F\in \mathfrak{T}(\Q)$ which is called the Frobenius at infinity.  For any $\sigma\in\Gal(\overline\Q/\Q)$, we let ${}^{\sigma}\mathfrak{S}=\pi^{-1}(\sigma)$. There is a cocharacter $\mu_{\can}:\Gm{\C}\to\mathfrak{S}_{\C}$, which in Tannakian terms gives rise to the Hodge cocharacter of the Hodge structures on $H_B(M)\otimes_{\Q}\C$.

Let $G$ be any algebraic group over $\Q$ and $\rho:\mathfrak{S}\to G^{\ad}$ be a homomorphism, inducing an action of $\mathfrak{S}$ on $G$ by group automorphisms (conjugation). Let ${}^{\sigma,\rho}G={}^{\sigma}\mathfrak{S}\times_{\mathfrak{S},\rho}G$ be group obtained by twisting $G$ by the torsor ${}^{\sigma}\mathfrak{S}$. Thus, ${}^{\sigma,\rho}G$ is the fpqc sheaf associated to the presheaf sending a $\Q$-algebra $R$ to the group ${}^{\sigma}\mathfrak{S}(R)\times_{\mathfrak{S}(R),\rho}G(R)$, which is the quotient of ${}^{\sigma}\mathfrak{S}(R)\times G(R)$ by the right action $(s,g)s_{1}=(ss_{1},s_{1}^{-1}g)$ of $\mathfrak{S}(R)$. We now specialize to the case $\sigma=c$. In this case, $^{c}\mathfrak{S}$ is already trivialized over $\Q$, a rational element being $sp(c)$. In particular, the map $sp(c)_{R}.g\mapsto g$ (for $g\in G(R)$) defines a group isomorphism between the above presheaf and $G$, and a fortiori between $^{c,\rho}G$ and $G$. If $H\subset G$ is a subgroup on which $\mathfrak{S}$ acts trivially, then $^{c,\rho}H$ is canonically isomorphic to $H$ (this is true for any $\sigma$), and the identifications are compatible. The isomorphism between $G(\A_f)$ and ${}^{c,\rho}G(\A_f)$, denoted by $g\mapsto{}^cg$ in \cite[II.4]{milne} becomes the identity under our identification, and similarly for the isomorphism $g\mapsto{}^cg$ between $G_{\C}$ and $^{c,\rho}G_{\C}$ defined in \cite[III.1]{milne} (note that the element $z_{\infty}(c)$ defined in {\em op. cit.} is equal to $sp(c)_{\C}$).

Suppose that $(G,X)$ is a Shimura datum as before, and $(T,x)$ is a special pair. The map $\mu_{x}$ factors through $T_{\C}$, and there exists a unique homomorphism $\rho_{x}^{\ad}:\mathfrak{S}\to G^{\ad}$ such that $(\rho_{x}^{\ad})_{\C}\circ\mu_{\can}=\mu_x^{\ad}$. For $\sigma\in\Aut(\C)$, the group ${}^{\sigma,x}G$ is defined to be ${}^{\sigma,\rho_x^{\ad}}G$ in the notation of the previous paragraph, where we take the restriction of $\sigma$ to $\overline\Q$. Since the cocharacter $\sigma(\mu_x)$ of $T={}^{\sigma,\rho_{x}}T$ commutes with its complex conjugate, it is the Hodge cocharacter associated to a map $\mathbb{S}\to{}^{\sigma,x}G_{\R}$ which we denote by ${}^{\sigma}x$, and ${}^{\sigma,x}X$ is defined to be its ${}^{\sigma,x}G(\R)$-conjugacy class. In particular, ${}^{c,x}G$ can be identified with $G$, and ${}^{c,x}X$ is naturally identified with the $G(\R)$-conjugacy class of ${}^{c}x$. Note that $^cx=x\circ c$, and $h\mapsto h\circ c$ defines an antiholomorphic isomorphism between $X$ and $^{c,x}X$. This doesn't  depend on $x$, and from now on we denote by $\overline X=\{h\circ c|h\in X\}$, and so the pair $(^{c,x}G,^{c,x}X)$ becomes naturally identified with the pair $(G,\overline X)$. The isomorphism $\varphi_{c,x}$ becomes, under this identification, an isomorphism $\varphi:\Sh_K(G,X)_E\times_{E,c}E\to\Sh_K(G,\overline X)_E$; in complex points, it defines an antiholomorphic isomorphism between $\Sh_K(G,X)(\C)$ and $\Sh_{K}(G,\overline X)(\C)$, which we denote by $\phi$. For $[h,g]\in\Sh_K(G,X)(\C)$, we have that $\phi([h,g])=[h\circ c,g]\in\Sh_K(G,\overline X)(\C)$.

For example, suppose that $E\subset\R$. Then there is an antiholomorphic involution on $\Sh_K(G,X)(\C)$ defined by complex conjugation acting on $\C$. It follows from the theory of canonical models that this involution takes the form $[h,g]\mapsto[\eta(h),g]$, where $\eta:X\to X$ is an antiholomorphic involution of the form $\eta(g.x)=(gn).x$ for some $n\in N(\R)$ (here $N$ is the normalizer in $G$ of $T$). In fact, the output is that there exists $n\in N(\R)$ such that ${}^cx=n.x$, and thus $\overline X=X$; then the map $\eta$ becomes what we called $\phi$, that is, $\eta(h)=h\circ c$ for any $h\in X$. 

\subsection{Involutions of Shimura data and descent} Fix a Shimura datum $(G,X)$, with reflex field $E$. For an involution $\theta:G\to G$, let $\theta(X)$ be the $G(\R)$-conjugacy class $\{\theta(h)\mid h\in X\}$, where $\theta(h)=\theta_{\R}\circ h$. Since we want to consider involutions $\theta$ that send $X$ to $\overline{X} \neq X$, from now on, we will focus on the case where $E$ is a CM field (if $E$ is totally real, the identity in $G$ takes $X$ to $\overline{X}$). Let $E^{+}\subset E$ be the maximal totally real subfield, and let $\iota\in\Gal(E/E^{+})$ be the non-trivial automorphism, i.e. the restriction of complex conjugation $c$ to $E$.

Suppose that $\theta$ is an involution of $G$ such that $\theta(X)=\overline X$. For a compact open subgroup $K\subset G(\A_f)$, denote by $^\theta K=\theta(K)\subset G(\A_f)$. Then $\theta$ induces an isomorphism of algebraic varieties $\Sh(\theta):\Sh_K(G,X)_E\to\Sh_{^\theta K}(G,\overline X)_E$. In complex points, this takes $[h,g]$ to $[\theta(h),\theta(g)]$. Suppose that $^\theta K=K$. Then $\Sh(\theta)^{-1}\circ\varphi$ defines an isomorphism $\psi:\iota(\Sh_{K}(G,X)_{E})=\Sh_K(G,X)_E\times_{E,\iota}E\to\Sh_K(G,X)_E$.

Let $V$ be an arbitrary scheme over $E$. Recall that an $E/E^{+}$-descent datum is a pair of isomorphisms $\psi_{1}:1(V)=V\times_{E,1}E\to V$ and $\psi_{\iota}:\iota V=V\times_{E,\iota}E\to V$ of schemes over $E$ satisfying the cocycle condition
\[ \psi_{\sigma}\circ\sigma(\psi_{\tau})=\psi_{\sigma\tau} \]
for all $\sigma,\tau\in\Gal(E/E^{+})$, using the natural identification $\sigma(\tau(V))=(\sigma\tau)V$. Then necessarily $\psi_{1}$ is the first projection $1(V)\to V$ and thus to give a descent datum amounts to give an isomorphism $\psi=\psi_{\iota}:\iota(V)\to V$ such that $\psi\circ\iota(\psi):\iota(\iota(V))\to V$ is equal to the identity map, when identifying $\iota(\iota(V))=V$. By definition, such a descent datum is effective if there exists a scheme $V_{0}$ over $E^{+}$ and an isomorphism $m:V\to V_{0,E}=V_{0}\times_{E^{+}}E$ such that $m\circ\psi=\iota(\psi)$, after identifying $\iota(V_{0,E})=V_{0,E}$. The descent criterion (\cite{Weil}) tells us that if $V$ is a quasi-projective algebraic variety, then a descent datum for $V$ is effective.

\begin{thm}\label{thmdescent} The map $\psi:\iota(\Sh_{K}(G,X)_{E})\to\Sh_{K}(G,X)_{E}$ obtained as above from an involution $\theta:G\to G$ such that $\theta(X)=\overline X$ and ${}^{\theta}K=K$  is an effective $E/E^{+}$-descent datum on the Shimura variety $\Sh_{K}(G,X)_{E}$. Hence, there exists a quasi-projective, smooth, algebraic variety $\Sh_K(G,X)_{E^+}$ over $E^+$, and an isomorphism $m:\Sh_K(G,X)_E\to\Sh_K(G,X)_{E^+}\times_{E^+}E$ such that $m\circ\psi=cm$. 
\begin{proof} Let $V=\Sh_{K}(G,X)_{E}$, $\overline V=\Sh_{K}(G,\overline X)_{E}$, and let $n:V\to\iota(\iota(V))$ be the natural isomorphism. We need to check that $\psi\circ\iota(\psi)\circ n=\id_{V}$, and for this it is enough to see that both morphisms are equal on the set of complex points $V(\C)$. Let $c_{V}:V(\C)\to(\iota V)(\C)$ be the bijection that sends $x:\Spec(\C)\to V$ to $p_{\iota,V}^{-1}\circ x\circ\Spec(c)$, where $p_{\iota,V}:\iota V\to V$ is the first projection, and define $c_{\iota V}:(\iota V)(\C)\to(\iota(\iota V))(\C)$ similarly. Then we have that $n(\C)=c_{\iota V}\circ c_{V}$, $\iota(\psi)(\C)=c_{V}\circ\psi(\C)\circ c_{\iota V}^{-1}$, and $\psi$ satisfies that $\psi(\C)\circ c_{V}=\Sh(\theta)^{-1}(\C)\circ\phi$. Recall that $\phi:V(\C)\to\overline V(\C)$ sends $[h,g]$ to $[\bar h,g]$. Putting all this together, we get that 
\[ (\psi\circ\iota(\psi)\circ n)(\C)=\Sh(\theta)^{-1}(\C)\circ\phi\circ\Sh(\theta)^{-1}\circ\phi, \]
and thus 
\[ (\psi\circ\iota(\psi)\circ n)(\C)([h,g])=[\theta^{-1}\left(\overline{\theta^{-1}(\overline h)}\right),\theta^{-2}(g)].\] But for any $y\in\overline X$, $\theta^{-1}(y)=\theta_{\R}^{-1}\circ y$, and so 
\[ \theta^{-1}\left(\overline{\theta^{-1}(\overline h)}\right)=\theta^{-1}\left(\overline{\theta_{\R}^{-1}\circ \overline h}\right)=\theta^{-1}\left(\overline{\theta_{\R}^{-1}\circ h\circ c}\right)=\theta^{-1}\left(\theta_{\R}^{-1}\circ h\right)=\theta_{\R}^{-2}\circ h=\theta^{-2}(h),\] and thus $(\psi\circ\iota(\psi)\circ n)(\C)([h,g])=[h,g]$, using the fact that $\theta^{2}=\id$. Finally, since $\Sh_K(G,X)_E$ is quasi-projective, the descent datum just constructed is effective.
\end{proof}
\end{thm}

\begin{rem}
The model of Theorem \ref{thmdescent} depends on the descent datum, which in turns depends on the particular involution $\theta$.
\end{rem}

We note that, by the nature of the descent datum, Hecke operators do not descend to the model $\Sh_{K}(G,X)_{E^{+}}$. Given $q\in G(\A_f)$, the Hecke operator $T_{q}$ is a morphism of algebraic varieties $T_q:\Sh_K(G,X)_E\to\Sh_{q^{-1}Kq}(G,X)_E$, which in complex points is given by $T_q([h,g])=[h,gq]$. Then $T_{\theta(q)}\circ \Sh(\theta)=\Sh(\theta)\circ T_q:\Sh_K(G,X)_E\to\Sh_{\theta(q)^{-1}{}^\theta K\theta(q)}(G,\overline X)_E$. The Hecke operator $T_q$ descends to a map $\Sh_K(G,X)_{E^+}\to\Sh_{q^{-1}Kq}(G,X)_{E^+}$ if and only if $T_{\theta(q)}=T_q$.

In the following sections we will construct several examples of involutions $\theta$ as above, and explain a general framework for such constructions.

\section{Opposition involutions}\label{sec:opposition}
In this section we recall some basic facts about opposition involutions and prove a few results that will be needed in the forthcoming sections. For the basic facts regarding root data, see \cite{springercorvallis}.

\subsection{Root data}\label{subsec:rootdata} Let $\Psi=(X,\Phi,X^\vee,\Phi^\vee)$ be a root datum with $\Phi\neq\emptyset$. Let $Q$ be the subgroup of $X$ generated by $\Phi$, and $V=Q\otimes_{Z}\Q$. Let $W=W(\Phi)$ be the Weyl group of the root system $\Phi$ in $V$. This can be naturally identified with the Weyl group of $\Phi^\vee$ and with the subgroup of $\Aut_{\Z}(X)$ generated by the reflections $s_{\alpha}$ for $\alpha\in\Phi$. Choose a basis $\Delta$, and consider the associated based root datum $\Psi_0=(X,\Phi,\Delta,X^\vee,\Phi^\vee,\Delta^\vee)$.

There is an obvious notion of isomorphism of root data (resp. based root data) $\Psi\to\Psi'$ (resp. $\Psi_0\to\Psi_0'$). It amounts to giving a $\Z$-linear isomorphism $f:X\to X'$ such that $f(\Phi)=\Phi'$ and ${}^tf(f(\alpha)^\vee)=\alpha^\vee$ for all $\alpha\in\Phi$ (resp. and $f(\Delta)=\Delta'$). Here ${}^tf$ denotes the transpose with respect to the root data pairings. We denote by $\Aut(\Psi)$ (resp. $\Aut(\Psi_0)$) the group of automorphisms of $\Psi$ (resp. $\Psi_0$). Each $s_{\alpha}$ can be seen as an automorphism of $\Psi$, and thus there is a natural inclusion $W\subset\Aut(\Psi)$. We also denote by $-1\in\Aut(\Psi)$ the automorphism that sends $x\in X$ to $-x\in X$.

Assume from now on that $\Phi$ is reduced. If $\Delta$ is a basis, let $w_0$ be the longest element of $W$ with respect to it. Then $w_0(\Delta)=-\Delta$, and thus $-w_0=-1\circ w_0\in\Aut(\Psi_0)$. We call $\star=-w_0$ the {\em opposition involution} of $\Psi_0$ (since $w_0^2=1$ it is indeed an involution). We denote the action of $\star$ on elements $x$ (which can be characters of $T$, nodes of the Dynkin diagram, etc.) by $x\mapsto x^\star$.  When $\Phi=\emptyset$, in which case $\Psi$ is called toral, we directly define $\star=-1\in\Aut_{\Z}(X)$.

\begin{rem} An isogeny (in particular, an isomorphism) of based root data will commute with the corresponding opposition involutions. In particular, $\star$ is a central element of $\Aut\Psi_0$.
\end{rem}

\begin{rem}\label{remtrivial} Let $X_{0}\subset X$ denote the subgroup of $X$ orthogonal to $\Phi^{\vee}$. The root datum $\Psi$ is called semisimple when $X_{0}=0$. If this is not the case, then there exists a non-zero $x\in X_{0}$, which hence must be invariant under $W$. In particular, $x^{\star}=-x\neq x$, so $\star$ cannot be the identity map if the root datum is not semisimple. In the same vein, if the root datum is toral then $\star\neq 1$ unless $\Psi$ is trivial (that is, also semisimple).
\end{rem}

Suppose now that $k$ is an algebraically closed field of characteristic $0$, and let $G$ be a reductive group over $k$. Let $T\subset G$ be a maximal torus, and $\Psi=\Psi(G,T)$ be the associated root datum, so that $X=X^{*}(T)$. Let $B\supset T$ be a Borel subgroup, and let $\Psi_0=\Psi_0(G,T,B)$ be the corresponding based root datum. Let $\Aut(G)$ be the group of automorphisms of $G$, and $\Inn(G)\subset\Aut(G)$ be the subgroup of inner automorphisms (that is, defined by elements in $G(k)$). Thus, $\Inn(G)\simeq G^{\ad}(k)\simeq G(k)/Z(k)$, where $Z$ is the center of $G$. Then there is a split exact sequence
\begin{equation}\label{exseq} 1\to\Inn(G)\to\Aut(G)\to\Aut\Psi_0\to 1 \end{equation}
where, for $f\in\Aut(G)$, the third arrow sends $f$ to the automorphism of $\Psi_0$ induced by $f'\in\Aut(G,T,B)$, where $f'=\inte(g)\circ f$ for any element $g\in G(k)$ such that $\inte(g)f(B,T)=(B,T)$. We define an {\em opposition involution} of $G$ (with respect to $(B,T)$) to be any element $\theta\in\Aut(G)$ of order 1 or 2 that induces the opposition involution $\star$ in $\Aut\Psi_0$. Note that this definition does not require $\theta$ to preserve $T$ or $B$. If $\theta'$ is another such involution then $\theta'=\inte(g)\circ\theta$ for some $g\in G(k)$. If $\theta$ is an opposition involution for $(B,T)$ and $(B',T')$ is another Borel pair, then it is also an opposition involution for $(B',T')$. 
The exact sequence (\ref{exseq}) is split by the choice of a pinning. More precisely, let $\Delta\subset\Phi$ be the set of simple roots corresponding to $B$. For each $\alpha\in\Delta$, let $U_\alpha\in G$ be the usual unipotent subgroup (see \cite{springercorvallis}, 2.3), and let $u_{\alpha}\in U_{\alpha}$ be a non-trivial element. The  {\em pinning} is the datum $\{u_{\alpha}\}_{\alpha\in\Delta}$ with respect to $(B,T)$, and a splitting $\Aut\Psi_0\to\Aut(G)$ of (\ref{exseq}) associated to this pinning is given by an isomorphism $\Aut\Psi_0\simeq\Aut(G,T,B,\{u_{\alpha}\}_{\alpha\in\Delta})$; two such splittings differ by an automorphism $\inte(t)$ for some $t\in T(k)$. In particular, after choosing a pinning, we can take $\theta\in\Aut(G)$ to be the image of $\star$ under the splitting and this will be an opposition involution, which proves their existence. Note that we are actually showing that there are opposition involutions in $\Aut(G)$ which preserve $T$ and $B$ (and a fixed pinning).

Let $k$ be any field of characteristic $0$, and $\overline k$ be an algebraic closure of $k$. Let $\Gamma=\Aut(\overline k/k)$. Let $G$ be a reductive group over $k$, $T\subset G$ a maximal torus, and $B\supset T_{\overline k}$ a Borel subgroup of $G_{\overline k}$. Let $\Psi=\Psi(G_{\overline k},T_{\overline k})$ and $\Psi_0=\Psi_0(G_{\overline k},T_{\overline k},B)$. There is a natural action of $\Gamma$ on $X$, denoted by $\chi\mapsto\chi^\gamma$, where 
\[ \chi^\gamma(t)=\gamma\left(\chi(\gamma^{-1}(t))\right)\]
for $\gamma\in\Gamma$, $t\in T(\overline k)$. We call it the {\em usual} action of $\Gamma$ on $X$. It defines an action of $\Gamma$ on $\Psi$. Let $\gamma\in\Gamma$. Then we define a second action $\mu_G(\gamma)$ on $X$, the $*$-action, given by $\mu_G(\gamma)(\chi)(t)=\chi^\gamma(n^{-1}tn)$ for $t\in T(\overline k)$, where $n\in G(\overline k)$ is an element such that $\inte(n)$ sends the Borel pair $(\gamma(B),\gamma(T_{\overline k}))$ to $(B,T_{\overline k})$. For example, if $B$ is a Borel defined over $k$, then we can take $n=1$ and the $*$-action is just the usual action $\chi\mapsto\chi^\gamma$. Going back to the general case, this gives a morphism $\mu_G:\Gamma\to\Aut\Psi_0$, and it induces an action of $\Gamma$ on $\Aut\Psi_{0}$ by taking $\rho\mapsto\mu_{G}(\gamma)\circ\rho\circ\mu_{G}(\gamma)^{-1}$ for $\rho\in\Aut\Psi_{0}$. There is also an action of $\Gamma$ on $\Aut(G_{\overline k})$ given by $\gamma.f=(1_G\times_{\Spec(k)}\Spec(\gamma^{-1}))\circ f\circ(1_G\times_{\Spec(k)}\Spec(\gamma))$, which on $G(\overline k)$-points is simply $g\mapsto\gamma(f(\gamma^{-1}(g)))$. It preserves the subgroup $\Inn(G_{\overline k})=G(\overline k)/Z(\overline k)$, where it acts as usual. The exact sequence (\ref{exseq}) becomes
\begin{equation}\label{exseqk} 1\to \Inn(G_{\overline k})\to\Aut(G_{\overline k})\to\Aut\Psi_0\to 1\end{equation}
and is $\Gamma$-equivariant. We define an opposition involution of $G$ to be an automorphism $\theta\in\Aut(G)$ of order 1 or 2 such that $\theta_{\overline k}$ is an opposition involution on $G_{\overline k}$. 

There may not be a $\Gamma$-equivariant splitting of (\ref{exseqk}), so it may not always be possible to construct in this way an opposition involution of $G$. However, if $G$ is quasi-split and $B$ is a Borel subgroup defined over $k$, it can be shown (see \cite{sga3},XXIV.3.10) that there exists a $\Gamma$-equivariant splitting. Since $\star\in\Aut\Psi_{0}$ is central, it commutes with $\mu_G(\gamma)$ for any $\gamma\in\Gamma$, and thus it is a $\Gamma$-invariant element in the last group of (\ref{exseqk}).  Thus, for quasi-split reductive groups over $k$, there always exist opposition involutions on $G$ over $k$,  but the condition of $G$ being quasi-split is far from necessary. There are many non-quasi-split cases where the opposition involution is trivial (see below), and so obviously defined over $k$. There are many non-trivial examples as well, as we will see later.

\begin{rem} If $G=T$ is a torus, then there exists one and only one opposition involution $\theta\in\Aut(G)$, namely $\theta=i_{G}$.
\end{rem}

\begin{lemma}\label{oppcenter} If $\theta$ is an opposition involution of $G$, then $\theta_{Z}:Z\to Z$ is equal to $i_{Z}$.
\begin{proof} It is enough to see that both maps induce the same map on $X^{*}(Z)$, that is, that $\theta_{Z}^{*}:X^{*}(Z)\to X^{*}(Z)$ is multiplication by $-1$, and thus we can assume that $k=\overline k$. Let $(B,T)$ be a Borel pair. Then $Z\subset T$. Let $\chi\in X^{*}(Z)$. Then there exists $\mu\in X^{*}(T)$ such that $\mu|_{Z}=\chi$. We claim that $\theta_{Z}^{*}(\chi)=(\mu^{\star})|_{Z}$. Indeed, for $z\in Z(k)$, $\theta_{Z}^{*}(\chi)(z)=\chi(\theta(z))$, whereas $(\mu^{\star})|_{Z}(z)=\Psi_{0}(\theta)(\mu)(z)=\mu((\inte(g)\circ\theta)(z))=\mu(\theta(z))$ (where $g\in G(k)$ sends $\theta(B,T)$ to $(B,T)$), which shows that $\theta_{Z}^{*}=(\mu^{\star})|_{Z}$. 

On the other hand, if $n_{0}\in N_{G}(T)(k)$ represents $w_{0}\in W=N_{G}(T)(k)/T(k)$, then for $z\in Z(k)$, $\mu^{\star}(z)=\mu(n_{0}^{-1}z^{-1}n_{0})=\mu(z^{-1})=\mu^{-1}(z)$ because $z\in Z(k)$. Thus, $\theta_{Z}^{*}(\chi)=-\chi$, as desired, where we have switched back to the additive notation for the group $X^{*}(Z)$.\end{proof}
\end{lemma}

\begin{rem} The last lemma shows in particular that if the identity is an opposition involution, then $Z$ is killed by $2$. Then $Z^{0}$ must be trivial, that is, $G$ must be semisimple (see also Remark \ref{remtrivial}).
\end{rem}

\subsection{Dynkin diagrams and special nodes}\label{subsec:dynkin}
Let $\Psi_{0}$ be a based root datum with $\Phi\neq\emptyset$ and reduced, and let $\mathcal{D}$ be its Dynkin diagram. Then the opposition involution $\star$ acts on $\mathcal{D}$. We include for reference the list of connected Dynkin diagrams and their opposition involutions; see \cite{bourbaki} for notation of nodes and more details. We also list the special nodes of each diagram (see \cite{delignecorvallis}, 1.2.5, for the definition of special node). Also, note that if $\Psi$ is semisimple, then $\star$ is trivial on $\Psi_{0}$ if and only if it is trivial on $\mathcal{D}$. For a Shimura datum $(G,X)$, the only factors of $G^{\ad}$ that contribute to a CM reflex field are the ones of type $A_{l}$ ($l\geq 2$), $D_{l}$ ($l\geq 5$ odd) or $E_{6}$. This follows from the list below and Proposition 2.3.6 of \cite{delignecorvallis}.
\begin{itemize}
\item $\mathcal{D}=A_l$ ($l\geq 1$)

 $\alpha_i^{\star}=\alpha_{l+1-i}$ (so $\star$ it's trivial if $l=1$)
 
 All nodes $\alpha_{i}$ are special

\item $\mathcal{D}=B_{l}$ ($l\geq 2$) or $C_{l}$ ($l\geq 3$)

  $\star$ is trivial
  
  There is only one special node: $\alpha_{1}$ in the $B_{l}$ case, and $\alpha_{l}$ in the $C_{l}$ case

\item $\mathcal{D}=D_{l}$ ($l\geq 4$)

 If $l$ is even, $\star$ is trivial
 
 If $l$ is odd, $\alpha_i^{\star}=\alpha_{i}$ for $i<l-1$, and $\alpha_{l-1}^{\star}=\alpha_{l}$
 
 The special nodes are $\alpha_{1}$, $\alpha_{l-1}$ and $\alpha_{l}$
 
\item $\mathcal{D}=E_6$

$\alpha_1^{\star}=\alpha_6$, $\alpha_2^{\star}=\alpha_{2}$, $\alpha_3^{\star}=\alpha_5$, $\alpha_4^{\star}=\alpha_4$

The special nodes are $\alpha_{1}$ and $\alpha_{6}$

\item $\mathcal{D}=E_7$, $E_8$, $F_4$ or $G_2$

$\star$ is trivial

Only $E_{7}$ has a special node, which is $\alpha_{7}$
\end{itemize}

\subsection{Multiplicative groups of CM type}\label{subsec:multgroups}

From now on let $k=\Q$ and $\Gamma=\Gal(\overline\Q/\Q)$. Let $T_{1}$ and $T_{2}$ be algebraic groups over $\Q$ of multiplicative type, not necessarily connected. Then there is a natural bijection $\Hom(T_{1},T_{2})\simeq\Hom_{\Gamma}(X_{2},X_{1})$, where $\Aut_{\Gamma}$ means $\Gamma$-equivariant morphisms for the natural Galois structures on $X_{i}=X^{*}(T_{i})$. In particular, for $T$ over $\Q$ of multiplicative type, there is a natural isomorphism $\Aut(T)\simeq\Aut_{\Gamma}(X)$, with $X=X^{*}(T)$. We let $c_{T}^{*}:X\to X$ be the map $c_{T}^{*}(\chi)=\chi^{c}$. We say $T$ splits over an extension $K\subset\overline\Q$ of $\Q$ if $\Aut(\overline\Q/K)$ acts trivially on $X^{*}(T)$. 

\begin{lemma} If $T$ is a group of multiplicative type that splits over a CM field, then $c_{T}^{*}\in\Aut_{\Gamma}(X)$.
\begin{proof} Suppose that $T$ splits over $K\subset\overline\Q$, a CM field. Let $\chi\in X$. Then $\chi^{\gamma}=\chi$ for any $\gamma\in\Aut(\overline\Q/K)$, and thus $\chi^{\gamma_{1}}=\chi^{\gamma_{2}}$ if $\gamma_{1},\gamma_{2}\in\Gamma$ have the same restriction to $K$. For any $\gamma\in\Gamma$, $\gamma c$ and $c\gamma$ have the same restriction to $K$, and so $c_{T}^{*}(\chi^{\gamma})=(\chi^{\gamma})^{c}=\chi^{c\gamma}=\chi^{\gamma c}=(c_{T}^{*}(\chi))^{\gamma}$.
\end{proof} 
\end{lemma}

Under the assumptions of the last lemma, we let $c_{T}:T\to T$ the unique involution inducing $c_{T}^{*}$ on $X$. If $T_{1}$ and $T_{2}$ are groups of multiplicative type which are split over a CM field, and $f:T_{1}\to T_{2}$ is a morphism, then $f\circ c_{T_{1}}=c_{T_{2}}\circ f$, because both maps induce the same morphism $X_{2}\to X_{1}$.

Suppose now that $T$ is a group of multiplicative type over $\R$. Using the same procedure, there exists a unique involution $c_{T}:T\to T$ inducing complex conjugation on characters. If $T$ is defined over $\Q$ and split over a CM field, these definitions are compatible with base change from $\Q$ to $\R$. 

\begin{ex} For $T=\mathbb{S}$ over $\R$, the map $c_{\mathbb{S}}$ is given by $c_{\mathbb{S}}(z\otimes a)=\overline z\otimes a$ for an $\R$-algebra $A$ and $z\otimes a\in(\C\otimes_{\R}A)^{\times}$.
\end{ex}

\begin{rem}\label{anistorus} If $T$ is an anisotropic $\R$-torus (that is, if $T(\R)$ is compact), then it's easy to see that $\chi^{c}=-\chi$ for any $\chi\in X$ and thus $c_{T}=i_{T}$ is the opposition involution on $T$.
\end{rem}

\subsection{Involutions taking $X$ to $\overline X$}\label{subsec:invoXtobarX}

Let $(G,X)$ be a Shimura datum. Recall that we are assuming that $Z^{0}$ splits over a CM field, and hence we have the conjugation involution $c_{Z^{0}}:Z^{0}\to Z^{0}$.

\begin{rem} Let $x\in X$. From the fact that $\inte(x(i)):G^{\ad}_{\R}\to G^{\ad}_{\R}$ is a Cartan involution, it follows that $G^{\ad}_{\R}$ is an inner form of an anisotropic group $H$ over $\R$ (that is, $H(\R)$ is compact). A similar statement holds for $G_{\R}^{\der}$ (the element $x(i)$ may not belong to $G^{\der}(\R)$, however over $\C$, $\inte(x(i))$ can be replaced by $\inte(x(i)')$ for some $x(i)'\in(T\cap G^{\der})(\C)$). The next lemma is well known.
\end{rem}

\begin{lemma}\label{anisotr} Let $G$ be a reductive group over $\R$ and assume that it is an inner form of a group $H$ over $\R$ which is anisotropic. Assume furthermore that $T\subset G$ is a maximal torus, and the inner automorphism of $G_{\C}$ defining a cocycle for $H$ is given by $\inte(t_{0})$ for some $t_{0}\in T(\C)$. Then the following hold.
\begin{enumerate}[(i)]
\item $c_{T}=i_{T}$.
\item For a Borel subgroup $B\supset T_{\C}$, the opposition involution acting on $\Psi_{0}(G,T,B)$ is given by the $*$-action of $c$.
\item The subgroup $c(B)\subset G_{\C}$ is the opposite Borel subgroup of $B$, that is, $c(B)\cap B=T_{\C}$.
\end{enumerate}
\begin{proof} Notice that $T(\R)$ is compact, so (i) follows from Remark \ref{anistorus}. By hypothesis, we can choose an isomorphism $\phi:G_{\C}\to H_{\C}$ such that $f:G_{\C}\to G_{\C}$ defined by $f(g)=\phi^{-1}\overline{\left(\phi(\overline g)\right)}$ is an inner automorphism of the form $\inte(t_{0})$, with $t_{0}\in T(\C)$. Then there exists a maximal torus $T_{H}\subset H$ such that $T_{H,\C}=\phi(T_{\C})$, and we let $B_{H}=\phi(B)$. Let $r=\Psi_{0}(\phi):\Psi_{0}(H,T_{H},B_{H})\to\Psi_{0}(G,T,B)$ be the induced isomorphism. It is $\Gal(\C/\R)$-equivariant for the $*$-actions, as follows from the fact that the forms are inner, and it commutes with $\star$, so it is enough to prove part (ii) when $G$ itself is anisotropic, which is well known. For (iii), the fact that $f$ preserves $T_{\C}$ and $B$ again allows us to reduce to the case of $G$ anisotropic, in which case the statement is well known.
\end{proof}
\end{lemma}

\begin{rem} In the last lemma, if the group if quasi-split and $B$ is a Borel subgroup defined over $\R$, the inner automorphism will not usually belong to $T(\C)$, otherwise we would have $B=T$. There are quasi-split semisimple groups with $B\neq T$ which are inner forms of anisotropic groups, for example $SU(n,n)$. In this case, the Cartan involution coming from a certain Shimura datum and special pair will preserve the maximal torus and a Borel subgroup containing it, but not a rational Borel subgroup.
\end{rem}

\begin{rem} Suppose that $(G,X)$ is a Shimura datum, and let $(T,x)$ be a special pair. Then $G_{\R}^{\der}$ satisfies all the hypotheses of the previous lemma. Here the inner automorphism defining the cocycle is $\inte(x(i)')$ as before. Alternatively, we can work with the adjoint group $G_{\R}^{\ad}$ and $x(i)$.
\end{rem}

\begin{rem}\label{remcinv} Suppose that $\theta:G\to G$ is an involution such that there exists a special pair $(T,x)$ with the property that $\theta$ preserves $T$ and induces $c_{T_{\R}}$ on $T_{\R}$. Then $\theta_{\R}(x)=c_{T_{\R}}\circ x=x\circ c_{\mathbb{S}}=\overline{x}$, and thus $\theta(X)=\overline X$.
\end{rem}

\begin{lemma}\label{conjborel} Let $G$ be a reductive group over $\R$, and $T\subset G$ a maximal torus. If $\theta:G\to G$ is an involution such that $\theta(T)\subset T$ and $\theta|_{T}=c_{T}$, then $\theta(B)=c(B)\subset G_{\C}$ for any Borel subgroup $B\supset T_{\C}$.
\begin{proof} Let $R\subset X=X^{*}(T)$ denote the set of roots of $(G_{\C},T_{\C})$. Let $R^{+}$ denote the set of positive roots with respect to $B$. Then $\theta(B)$ is the Borel subgroup whose Lie algebra is $\Lie(T_{\C})\oplus\bigoplus_{\alpha\in R^{+}}\Lie(G_{\C})_{\alpha\circ\theta}$. Since $\alpha\circ\theta=\alpha^{c}$, it follows that this is the Lie algebra of $c(B)$, and since both $\theta(B)$ and $c(B)$ are connected, this proves the lemma.
\end{proof}
\end{lemma}

The construction of involutions taking $X$ to $\overline X$ that we will perform will be based on involutions $\theta$ which will roughly be as in Remark \ref{remcinv}. By the following proposition, we need to look for opposition involutions on semisimple groups. 

\begin{prop} Let $(G,X)$ be a Shimura datum and let $\theta:G\to G$ be an involution of $G$, such that there exists a special pair $(T,x)$ with the property that $\theta$ preserves $T$ and induces $c_{T_{\R}}$ on $T_{\R}$. Then $\theta^{\der}:G^{\der}\to G^{\der}$ is an opposition involution, and $\theta_{0}=\theta|_{Z^{0}}:Z^{0}\to Z^{0}$ is equal to $c_{Z^{0}}$.
\begin{proof} Suppose that $\theta$ is an involution with $(T,x)$ as in the statement. To see that $\theta_{0}=c_{Z^{0}}$ it's enough to see that $\theta_{\R,0}=c_{Z^{0}_{\R}}$. Since $Z^{0}_{\R}\subset T_{\R}$ and $\theta_{\R}|_{T_{\R}}=c_{T_{\R}}$ it follows that $\theta_{0,\R}=c_{Z^{0}_{\R}}$. Let $T'=T\cap G^{\der}$, let $B\subset G_{\C}$ be a Borel subgroup containing $T_{\C}$, and $B'=B\cap G^{\der}_{\C}\supset T'_{\C}$. Let $\Psi_{0}'=\Psi_{0}(G^{\der},T',B')$ and let $r=\Psi_{0}(\theta^{\der}):\Psi_{0}'\to\Psi_{0}'$ be the induced isomorphism. It is given by $r(\chi)=\chi\circ\inte(q)\circ\theta^{\der}|_{T'_{\C}}$ for $\chi\in X'=X^{*}(T')$, where $q\in G^{\der}(\C)$ is such that $\inte(q)\theta^{\der}(T'_{\C},B')=(T'_{\C},B')$. On the other hand, by Lemma \ref{anisotr}, (ii), $\star:\Psi_{0}'\to\Psi_{0}'$ is given by $\chi^{\star}=\chi^{c}\circ\inte(a^{-1})$, where $a\in G^{\der}(\C)$ is such that $\inte(a)c(T'_{\C},B')=(T'_{\C},B')$. By Lemma \ref{conjborel}, we can take $a=q$. Finally, the hypothesis that $\theta^{\der}|_{T'}=c_{T'}$ implies that $\chi^{\star}=\chi\circ\theta\circ\inte(q^{-1})$. Thus, to see that $r(\chi)=\chi^{\star}$, it is enough to see that $\theta^{\der}\circ\inte(q^{-1})$ and $\theta^{\der}\circ\inte(\varphi(q))$ induce the same automorphism of $T'_{\C}$, and this follows from the fact that both elements $\theta^{-1}(q)$ and $q^{-1}$ conjugate the Borel pair $(T'_{\C},B')$ to the same Borel pair.
\end{proof}
\end{prop}

The following proposition is a partial converse and the main result of this section. Since our construction will be explicit using the classification of semisimple groups, we need to work with either the derived group or the adjoint group. The idea is to construct an involution on $G$ taking $X$ to $\overline{X}$ by extending an opposition involution on $G^{\der}$. Ideally we would want the involution to be as in Remark \ref{remcinv}, but it is enough to consider a weaker hypothesis, as stated in the proposition. Recall the notation from Subsection \ref{subsec:notshimura}. Suppose that for each $i$, $S_{i}\subset H_{i}$ is a maximal torus, and let $\widetilde T_{i}=\Res_{F_{i}/\Q}S_{i}\subset\widetilde G_{i}$, $T_{i}\subset G_{i}$ its image in $G_{i}$, $T'\subset G^{\der}$ the image of their product, and $T=Z^{0}T'$. Note that $T_{\R}^{\ad}=T_{\R}^{'\ad}=\prod_{i,v}S_{i,v}^{\ad}$, where $S_{i,v}\subset H_{i,v}$ and $S_{i,v}^{\ad}$ is its image in $H_{i,v}^{\ad}$.

\begin{prop}\label{mainprop} 

Suppose that $\theta_{i}:H_{i}\to H_{i}$ is an opposition involution for each $i$. Suppose moreover that $\theta_{i}(S_{i})=S_{i}$ and $\theta_{i,v}^{\ad}|_{S_{i,v}^{\ad}}=c_{S_{i,v}^{\ad}}$ for every $i$ and $v\in I_{i,nc}$. Finally, assume that there exists $x\in X$ such that $x^{\ad}$ factors through $T_{\R}^{\ad}$. Then there exists an involution $\theta:G\to G$ such that $\theta(X)=\overline X$.

\begin{proof} For each $i$, the involution $\Res_{F_{i}/\Q}\theta_{i}$ defines an opposition involution of $\widetilde G_{i}$. Moreover, the kernel $K_{i}$ of the projection $\widetilde G_{i}\to G_{i}$ is contained in the center of $\widetilde G_{i}$. By Lemma \ref{oppcenter}, $\Res_{F_{i}/\Q}\theta_{i}$ induces $x\mapsto x^{-1}$ on the center. In particular, it preserves $K_{i}$ and induces an opposition involution on $G_{i}$. Similarly, the product of these involutions define an opposition involution $\theta':G^{\der}\to G^{\der}$. Let $q:Z^{0}\times G^{\der}\to G$ be the natural isogeny. We can look at the product involution $\theta'\times c_{Z^{0}}:G^{\der}\times Z^{0}\to G^{\der}\times Z^{0}$. We claim that this preserves $\ker(q)$, and thus it induces an involution on $G$. To show this, we can work with $\C$-points. The kernel consists of pairs $(g,z)$ such that $zg=1$, so we need to check that if $(g,z)$ is such a pair, then $\theta'(g)c_{Z^{0}}(z)=1$.  The element $g=z^{-1}$ belongs to $Z^{0}\cap G^{\der}\subset Z_{G^{\der}}$. The maps $c_{Z^{0}}:Z^{0}\to Z^{0}$ and $c_{Z_{G^{\der}}}:Z_{G^{\der}}\to Z_{G^{\der}}$ are equal on $Z^{0}\cap G^{\der}$, and by part (i) of Lemma \ref{anisotr}, $c_{Z_{G^{\der}}}=i_{Z_{G^{\der}}}$, so $c_{Z^{0}}(z)=z^{-1}$. On the other hand, by Lemma \ref{oppcenter}, $\theta'$ induces $i_{Z_{G^{\der}}}$ on $Z_{G^{\der}}$, and so $\theta'(g)=g^{-1}=z$. This proves that there exists a (unique) involution $\theta:G\to G$ such that $\theta^{\der}=\theta'$ and $\theta_{0}=c_{Z^{0}}$.

We also have that $\theta$ preserves $T$ and $\theta_{\R}^{\ad}=\prod_{i,v}\theta_{i,v}^{\ad}$. Now, we know that there exists $x\in X$ such that $x^{\ad}$ factors through $T_{\R}^{\ad}$. Let $y=\theta_{\R}(x)$. Then $y^{\ad}=\theta_{\R}^{\ad}(x^{\ad})$. For $v\in I_{i,nc}$, we have $\theta_{i,v}^{\ad}(x_{i,v})=\overline{x_{i,v}}$ because $\theta_{i,v}^{\ad}$ induces $c_{S_{i,v}^{\ad}}$ and $x_{i,v}$ factors through $S_{i,v}^{\ad}$. For $v\in I_{i,c}$, $x_{i,v}=1$. Thus, $y^{\ad}=\overline{x}^{\ad}$. Also, since $\theta_{\R}$ induces $c_{Z^{0}_{\R}}$ on $Z^{0}_{\R}$, and $q:Z^{0}\to G^{\ab}$ is an isogeny, it follows that $\theta_{\R}$ induces $c_{G^{\ab}_{\R}}$ on $G^{\ab}_{\R}$. From this it follows that $y$ and $\overline x$ have the same projections to $G_{\R}^{\ad}$ and to $G^{\ab}_{\R}$, and thus $y=\overline x$ (see for instance the proof of Proposition 5.7 of \cite{milneIntro}). Since $y=\theta_{\R}(x)$, this shows that $\theta(X)=\overline X$.
\end{proof}
\end{prop}

\section{Involutions on classical semisimple groups}\label{sec:semisimple}
In this section, we make use of several results regarding the classification of semisimple algebraic groups over totally real fields. For notation and terminology regarding algebras with involutions and their associated groups, we freely follow our main reference \cite{involutions}. We are only interested in the explicit classification of groups of type $A$ and $D$ in order to construct our desired involutions on certain Shimura varieties. Furthermore, not all the groups in the general classification appear in the theory of Shimura varieties, so we are only interested in classifying the groups $H_{i}$ (in the notation of Subsection \ref{subsec:notshimura}) of type $A_{l}$ ($l\geq 2$) or $D_{l}$ ($l\geq 4$ odd) that can occur. Furthermore, in accordance with the previous section, we are also interested in constructing, whenever possible, opposition involutions on these groups.

The following construction regarding quaternion algebras will be used often in the following. Suppose that $D$ is a quaternion division algebra over a number field $K$. Let $\lambda\in D^{\times}$ be a pure quaternion (that is, such that $\sigma(\lambda)=-\lambda$, where $\sigma:D\to D$ is the canonical involution), and choose another pure quaternion $\mu\in D^{\times}$ such that $\lambda\mu=-\mu\lambda$. Then $\{1,\lambda,\mu,\lambda\mu\}$ is a standard basis of $D$. If we let $L=K(\lambda)$, then $L$ is a maximal subfield of $D$ (a quadratic extension of $K$). We have an isomorphism of $L$-algebras $\phi:D\otimes_{K}L\to M_{2}(L)$ defined by 
\[ \phi(\lambda\otimes1)=\left(\begin{array}{cc}\lambda & 0 \\ 0 & -\lambda\end{array}\right)\]
and
\[ \phi(\mu\otimes1)=\left(\begin{array}{cc} 0 & \mu^{2} \\ 1 & 0\end{array}\right).\]
Then the isomorphism $\phi$ sends $L\otimes_{K}L$ to the subalgebra of diagonal matrices in $M_{2}(L)$.

Throughout this section, let $F$ be a totally real field and $H$ be an absolutely almost simple, simply connected algebraic group over $F$. We let $\mathcal{D}$ be the Dynkin diagram of $H_{\overline F}$ (where $\overline F$ is some algebraic closure of $F$). We let $I=\Hom(F,\R)$, $I_{c}=\{v\in I:H_{v}^{\ad}(\R)\text{ is compact}\}$ and $I_{nc}$ its complement in $I$.

\subsection{Groups of type $A_{l}$ ($l\geq 2$)}\label{subsec:groupsAl}
Suppose that $\mathcal{D}=A_{l}$ with $l\geq 2$. Then there exists a quadratic \'etale extension $K/F$ (so $K/F$ is a quadratic extension of fields, or $K=F\times F$), and a central simple algebra $B$ over $K$, of degree $l+1$, endowed with an involution $\tau:B\to B$ of the second kind (that is, inducing $\iota$ on $K$, where $\iota$ is the non-trivial automorphism of $K$ which fixes $F$) such that $H=\mathbf{SU}(B,\tau)$ (\cite{involutions}, Theorem 26.9). If $H$ is one of the $H_{i}$ as above, then $K$ must be a field. Indeed, if otherwise, then $H\simeq\mathbf{SL}_{1}(A)$ for some central simple algebra $A$ over $F$ of degree $l+1$. For each $v\in\Hom(F,\R)$, $A_{v}=A\otimes_{F,v}\R\simeq M_{l+1}(\R)$ or $A_{v}\simeq M_{(l+1)/2}(\mathbb{H})$. In both cases, it follows that $H_{v}$ is an inner form of $\SL_{l+1,\R}$, so the $*$-action of $c$ is trivial (a condition that doesn't depend on the Borel pair), and thus it cannot be the opposition involution because $l\geq 2$. From this and Lemma \ref{anisotr} it follows that $H$ cannot occur as one the factors $H_{i}$. Thus, we have proved that $K$ must be a field. Moreover, a similar argument implies that $K$ must be totally imaginary, that is, $K/F$ is a CM extension. The adjoint group $H^{\ad}$ is $\mathbf{PGU}(B,\tau)$.

We can then write $B=\End_{D}(V)$ for some central division algebra $D$ over $K$, endowed with an involution $J:D\to D$ of the second kind, whose action we denote by $d\mapsto d^{J}$, and a finite dimensional right $D$-vector space $V$. There is a non-degenerate hermitian form $h:V\times V\to D$ inducing the involution $\tau:B\to B$. The pair $(V,h)$ is called a hermitian space over $D$.

Suppose that $\theta:H\to H$ is an opposition involution. There is a natural isomorphism between $\Aut(H)$ and the group of $F$-algebra automorphisms of $B$ that commute with $\tau$ (\cite{involutions}, Theorem 26.9), and thus there exists such an automorphism $\gamma:B\to B$ of order $2$, inducing $\theta$. If $\gamma|_{K}$ is the identity of $K$, then $\gamma=\inte(b_{0})$ for some $b_{0}\in B^{\times}$ by the Skolem-Noether Theorem, and $b_{0}$ is moreover a similitude for $\tau$. The induced map $\theta:H\to H$ would thus be an inner automorphism, inducing the identity on the Dynkin diagram, but the opposition involution on $A_{l}$ is non-trivial for $l\geq 2$. Hence, $\gamma|_{K}$ must be $\iota$. Let $\overline{B}$ and $\overline{D}$ denote the $K$-algebras $B$ and $D$ with $\iota$-conjugate structure. Thus $\gamma:B\to\overline B$ is a $K$-algebra isomorphism. We let $\operatorname{Br}(K)$ be the Brauer group of $K$ and $[B]=[D]\in\operatorname{Br}(K)$ be the class of $B$ in it. Then $[D]=[B]=[\overline B]=[\overline D]$, which implies that there must exist a ring automorphism $\alpha:D\to D$ inducing $\iota$ on $K$.

\begin{prop}\label{propcentral} Let $D$ be a central division algebra over a CM extension $K/F$ of number fields, endowed with an involution $J:D\to D$ of the second kind. Then the following are equivalent:
\begin{enumerate}[(a)]
\item $D=K$ or $D$ is a quaternion divison algebra over $K$.
\item The order of $[D]\in\operatorname{Br}(K)$ is $1$ or $2$.
\item There exists a ring automorphism $\alpha:D\to D$ inducing $\iota$ on $K$.
\end{enumerate}
Moreover, in this case, $\alpha$ is unique up to composition with an inner automorphism of $D$. Furthermore, it can be chosen to have order $2$ and such that $\alpha J=J\alpha$ is either $1$ if $D=K$ or the canonical involution if $D$ is a quaternion division algebra.
\begin{proof} The fact that $(a)$ implies $(b)$ in the quaternion algebra case follows from the existence of the canonical involution on $D$, which gives an isomorphism $D\to D^{\operatorname{op}}$, so $[D]=[D]^{-1}$. To see that $(b)$ implies $(a)$, see 10.2.3 of \cite{scharlau}.

Now suppose that $(a)$ is true. If $D=K$, then take $\alpha=\iota$. If $D$ is a quaternion division algebra, let $\sigma:D\to D$ be its canonical involution, and take $\alpha=J\sigma=\sigma J$ (they commute because $J\sigma J$ is a symplectic involution of the first kind on $D$, and hence equal to $\sigma$).

Finally, suppose that $\alpha:D\to D$ is as in $(c)$. The involution $J:D\to D$ induces an isomorphism $D\to\overline{D}^{\operatorname{op}}$, where $\overline{D}$ is the conjugate algebra $\lambda\cdot d=\iota(\lambda)d$ for $\lambda\in K$. Similarly, $\alpha$ induces an isomorphism $D\to\overline{D}$, and thus in the end we have an isomorphism $D\to D^{\operatorname{op}}$, which implies that the order of $[D]$ is $1$ or $2$.

The uniqueness of $\alpha$ up to inner automorphism follows because if $\beta$ is another such automorphism, then $\alpha\beta^{-1}:D\to D$ is a $K$-linear automorphism and hence inner by the Skolem-Noether Theorem.
\end{proof}
\end{prop}

\begin{rem}\label{remcentral} Suppose that $D$ is a quaternion division algebra. Under the conditions of the previous proposition, there exists a unique quaternion algebra $D_{0}\subset D$ over $F$ such that $D=D_{0}\otimes_{F}K$ and $J=\sigma_{0}\otimes_{F}\iota$, where $\sigma_{0}$ is the canonical involution of $D_{0}$ (see \cite{involutions}, 2.22). Then the map $\alpha$ constructed in the proof is $\alpha=1_{D_{0}}\otimes_{F}\iota$. We define the {\em canonical conjugation} $\alpha:D\to D$ (attached to $J$ or $D_{0}$) to be $\alpha=1_{D_{0}}\otimes_{F}\iota$. If $D=K$, we also call $\alpha=\iota$ the canonical conjugation.
\end{rem}

Thus, we have shown that if there exists $\theta:H\to H$ an opposition involution, then either $D=K$ (and $J=\iota$) or $D$ is a quaternion division algebra (and $J=\sigma_{0}\otimes_{F}\iota$). Conversely, suppose that $D=K$ or $D$ is a quaternion division algebra. We will construct suitable opposition involutions under an additional assumption. 

\begin{rem} Suppose that $D=K$ (and $D_{0}=F$) or $D$ is a quaternion division algebra. Let $I_{s}\subset\Hom(F,\R)$ be the subset of places $v\in I=\Hom(F,\R)$ such that $D_{0,v}=D_{0}\otimes_{F,v}\R$ is split, and let $I_{ns}\subset I$ be its complement. We let $I_{c}\subset I$ be the subset of places $v$ such that $H_{v}^{\ad}(\R)$ is compact, and $I_{nc}$ its complement. The group $H_{\R}$ can be written as a product of special unitary groups $\prod_{v\in I}\operatorname{SU}(p_{v},q_{v})$, and the compact places are exactly the places where $p_{v}q_{v}=0$. 

\end{rem} 

\begin{dfn}\label{dfn:strhermitian} We say that a hermitian space $(V,h)$ over $D$ (where $D=K$ or a quaternion division algebra) is {\em strongly hermitian} if there exists an $h$-orthogonal $D$-basis $\beta=\{v_{1},\dots,v_{n}\}$ of $V$ such that $h(v_{i},v_{i})\in K^{\times}$ for all $i$; in the quaternion algebra case, we ask furthermore that $I_{ns}\subset I_{c}$.
\end{dfn}

\begin{rem} Note that a hermitian space over $D=K$ is always strongly hermitian. 
\end{rem}

Suppose that $(V,h)$ is strongly hermitian, and let $\beta$ be a basis as in the definition. Let $I:V\to V$ be the $\alpha$-semilinear isomorphism obtained by applying $\alpha$ to the coordinates of elements of $V$ with respect to the basis $\beta$ (this map is inspired in the constructions of \cite{taylor}). Then $h(I(x),I(y))=\alpha(h(x,y))$. Let $\theta:H\to H$ be given as $\theta_{A}(g)=I_{A}gI_{A}$ for an $F$-algebra $A$ and a $D\otimes_{F}A$-linear automorphism $g$ of $V\otimes_{F}A$. Let $L\subset D$ be a maximal subfield. More precisely, if $D=K$ then $L=K$, and if $D$ is a quaternion division algebra, take $L=K(\lambda)$, where $\lambda$ is a pure quaternion in $D_{0}$. Let $S=S_{L,\beta}$ be the subgroup of $H$ defined as follows. For an $F$-algebra $A$, 
\[ S(A)=\{h\in H(A)\subset\Aut_{D\otimes_{F}A}(V\otimes_{F}A):h(v_{i}\otimes1)=(v_{i}\otimes 1)\lambda_{i},\text{ for some }\lambda_{i}\in(L\otimes_{F}A)^{\times}\}.\]
This is a maximal torus in $H$.

\begin{prop} With the above hypotheses, the following statements are true.
\begin{enumerate}[(a)]
\item The involution $\theta:H\to H$ is an opposition involution.
\item We have $\theta(S)=S$ and for every $v\in I_{nc}$, $\theta_{v}:H_{v}\to H_{v}$ induces $c_{S_{v}}$ on $S_{v}$.
\end{enumerate}
In particular, $\theta^{\ad}_{v}:H^{\ad}_{v}\to H^{\ad}_{v}$ induces $c_{S^{\ad}_{v}}$ on $S^{\ad}_{v}$ for $v\in I_{nc}$.
\begin{proof}
For part (a), it suffices to see that $\theta_{L}:H_{L}\to H_{L}$ is an opposition involution. We can identify $H_{K}$ with $SL_{V/D}$, where $SL_{V/D}(A)$ consists, for a $K$-algebra $A$, of the $D\otimes_{K}A$-linear automorphisms of $V\otimes_{K}A$ with reduced norm $1$. Using the basis $\beta$, we can further identify $SL_{V/D}(A)\cong\SL_{n}(D\otimes_{K}A)$. Let $Q\in\GL_{n}(K)$ be the matrix of $h$ with respect to $\beta$. Then it's easy to see that $\theta_{A}:\SL_{n}(D\otimes_{K}A)\to\SL_{n}(D\otimes_{K}A)$ is explicitly given by the formula
\[ \theta_{A}(X)=Q^{-1}({}^{t}X^{-1})^{\sigma}Q,\]
where $\sigma:D\to D$ is the canonical involution of $D$ if $D$ is a quaternion division algebra, and $\sigma=\id$ if $D=K$. Note that $Q$ is a diagonal matrix in $\GL_{n}(K)$.

If $D=K$, we denote by $\phi:D\otimes_{K}L\to L$ the unique obvious isomorphism. If $D$ is a quaternion division algebra, we take $\phi:D\otimes_{K}L\to M_{2}(L)$ to be an isomorphism of $L$-algebras taking $L\otimes_{K}L$ to $D_{2}(L)$ as constructed above (we use for this the pure quaternion $\lambda\in D_{0}$ and another pure quaternion $\mu\in D_{0}$ such that $\lambda\mu=-\mu\lambda$). In particular, $\sigma$ preserves $L$. The identification $H_{K}(A)\cong\SL_{n}(D\otimes_{K}A)$ sends $S_{K}(A)$ to the subgroup of matrices in $\SL_{n}(D\otimes_{K}A)$ which are diagonal and have entries in $L\otimes_{K}A$. Since $\sigma$ preserves $L$, it follows that $\theta_{K}$ sends the torus $S_{K}$ to itself. Moreover, if we now extend scalars to $L$, the map $\phi$ provides an isomorphism
\begin{equation}\label{identHL} H_{L}\cong\SL_{nm,L},\end{equation}
where $S_{nm,L}$ is the usual group of $nm\times nm$-matrices of determinant $1$; furthermore, the torus $S_{L}$ maps to the torus of diagonal matrices in $\SL_{nm,L}$ (so $S$ in indeed is a maximal torus, as claimed). If $m=1$, then $\theta_{L}(X)=Q^{-1}{}^{t}X^{-1}Q$ for $X\in\SL_{n,L}$. Suppose that $m=2$. Write $Q=\diag(q_{1},\dots,q_{n})$, and let $\widetilde Q=\diag(q_{1},\dots,q_{n},q_{1},\dots,q_{n})\in\GL_{2n}(K)$. Write matrices $X\in\SL_{2n,L}$ as blocks
\[ X=\left(\begin{array}{cc}A & B \\ C & D\end{array}\right),\]
with $A,B,C,D$ of size $n\times n$. Then $\theta_{L}:\SL_{2n,L}\to\SL_{2n,L}$ is explicitly given as
\[ \theta_{L}(X)=\widetilde Q^{-1}\left(\begin{array}{cc}{}^{t}D & -{}^{t}B \\ -{}^{t}C & {}^{t}A\end{array}\right)\widetilde Q.\]
From this explicit expression of $\theta$ as an involution of $\SL_{nm,L}$, it's easy to see that it preserves the maximal torus $S_{L}$ of diagonal matrices and that it induces the opposition involution on the root datum. 

For part (b), fix $v\in I_{nc}$. We need to check that if $\chi\in X=X^{*}(S_{v})=\Hom(S_{v}\times_{\R}\C,\Gm{\C})$, then $\chi\circ\theta_{v,\C}=\chi^{c}$. To compute $\chi^{c}$, we need to compute how complex conjugation acts on $H_{v}(\C)$. Choose once and for all an extension $\tau:L\hookrightarrow\C$ of $v$ to $L$. Using the embedding $\tau$ and the isomorphism (\ref{identHL}), we can identify $H_{v}\times_{\R}\C=H_{L}\times_{L,\tau}\C\cong\SL_{nm,\C}$. Moreover, the action of $c$ on $H_{v}(\C)\cong\SL_{nm}(\C)$ is explicitly given as follows. Let $Q_{v}=\diag(v(q_{1}),\dots,v(q_{n}))\in\GL_{n}(\R)$ and $\widetilde Q_{v}=\diag(v(q_{1}),\dots,v(q_{n}),v(q_{1}),\dots,v(q_{n}))\in\GL_{2n}(\R)$. Let 
\[ \gamma=\left(\begin{array}{cc} 0 & I_{n} \\ -I_{n} & 0\end{array}\right).\]
If $m=1$ and $X\in\SL_{n}(\C)$, then $c(X)=Q_{v}^{-1}X^{*,-1}Q$. If $m=2$ and $X\in\SL_{2n}(\C)$, then $c(X)=Q_{v}^{-1}\gamma X^{*,-1}\gamma^{-1}Q_{v}$. The last case easily follows from (\ref{identHL}) and the fact that $D_{0,v}$ is split. We can identify $X^{*}(S_{v})$ in the standard way with $\Z^{nm}/L$, where $L=\{(k,k,\dots,k):k\in\Z\}$. It then follows easily from our calculations of the action of $c$ that if $\chi\in X^{*}(S_{v})$ is identified with the class of the tuple $(a_{1},\dots,a_{n})$ in the case $m=1$, respectively the class of the tuple $(a_{1},\dots,a_{n},b_{1},\dots,b_{n})$ in the case $m=2$, then $\chi^{c}$ is identified with $(-a_{1},\dots,-a_{n})$ or with $(-b_{1},\dots,-b_{n},-a_{1},\dots,-a_{n})$ respectively. This, together with our formulas for $\theta$, show that $\theta_{v}$ induces $c_{S_{v}}$ on $S_{v}$, which is what we wanted to prove. 
\end{proof}
\end{prop}

\begin{rem} When $D_{0,v}$ is not split, there is also an explicit formula for $c$ which involves a matrix $\gamma$ as above, but which is a diagonal matrix. So in this case $\theta_{v}$ does not induce $c_{S_{v}}$ on $S_{v}$. We only care about non-compact places, hence our assumption $I_{ns}\subset I_{c}$.
\end{rem}

\begin{rem}\label{rem:mapyA} Keep the assumptions and notation as above. For each $v\in I_{nc}$, we will construct a map $y:\mathbb{S}\to H_{v}^{\ad}$ satisfying Deligne's axioms (\cite{delignecorvallis}, 1.2.1) and factoring through $S_{v}^{\ad}$. Namely, fix $\tau:L\hookrightarrow\C$ an extension of $v$ to $L$, and let $w=\tau|_{K}$ (so $w=\tau$ when $D=K=L$). Let $D_{w}=D\otimes_{K,w}\C$ and $J_{w}:D_{w}\to D_{w}$ be defined by $J_{w}(d\otimes z)=J(d)\otimes\overline z$. The group $H_{v}(A)$ can be identified, using the basis $\beta$, with the group of matrices $X\in\GL_{n}(D_{w}\otimes_{\R}A)$ such that ${}^{t}X^{J_{w}}QX=Q$ and $\operatorname{Nrd}(X)=1$. If $m=1$, let $\phi_{\tau}:D_{w}\to\C$ be the unique isomorphism. If $m=2$, consider the $\C$-algebra isomorphism $\phi_{\tau}:D_{w}\to M_{2}(\C)$ given by
\[ \phi_{\tau}(\lambda\otimes_{K,w}1)=\left(\begin{array}{cc}\tau(\lambda) & 0 \\ 0 & -\tau(\lambda) \end{array}\right),\]
\[ \phi_{\tau}(\mu\otimes_{K,w}1)=\left(\begin{array}{cc}0 & v(\mu^{2}) \\ 1 & 0\end{array}\right).\]
As above, for any $\R$-algebra $A$, this induces an isomorphism $\GL_{n}(D_{w}\otimes_{\R}A)\cong\GL_{mn}(\C\otimes_{\R}A)$ taking the subgroup of diagonal matrices with entries in $L_{w}\otimes_{\R}A$ (where $L_{w}=L\otimes_{K,w}\C$) to the subgroup of diagonal matrices in $\GL_{mn}(\C\otimes_{\R}A)$. Moreover, the corresponding involution $X\mapsto{}^{t}X^{J_{w}}$ gets identified with $X\mapsto\gamma X^{*}\gamma^{-1}$, where if $m=1$, $\gamma=I_{n}$, and if $m=2$, $\gamma$ is the hermitian matrix defined by
\[ \gamma=\left(\begin{array}{cc}0 & iI_{n}\\ -iI_{n} & 0\end{array}\right)\]
if $v(\lambda^{2})>0$, and
\[ \gamma=\left(\begin{array}{cc}-v(\mu^{2})I_{n} & 0 \\ 0 & I_{n}\end{array}\right)\]
if $v(\lambda^{2})<0$ (note that in this case, we must have $v(\mu^{2})>0$). In this way, we can write
\[ H_{v}(A)\cong\left\{X\in\GL_{2n}(\C\otimes_{\R}A):(\gamma X^{*}\gamma^{-1})Q'X=Q',\quad\det(X)=1\right\},\]
where $Q'=Q_{v}$ if $m=1$ and $Q'=\widetilde Q_{v}=\diag(v(q_{1}),\dots,v(q_{n}),\dots,v(q_{1}),\dots,v(q_{n}))$ if $m=2$. Thus, we can identify $H_{v}$ with the special unitary group $SU(\gamma^{-1}Q')$ of the hermitian matrix $\gamma^{-1}Q'$, and the maximal torus $S_{v}$ is the torus of diagonal matrices. Note that $H_{v}^{\ad}$ is also the adjoint group of the similitude unitary group $\operatorname{GU}(\gamma^{-1}Q')$. We define $y':\mathbb{S}\to\operatorname{GU}(\gamma^{-1}Q')$ as follows. For an $\R$-algebra $A$ and $z\in\mathbb{S}(A)$, let
\[ y_{A}'(z)=\left(\begin{array}{cc}\diag(y_{A}'(z)_{1},\dots,y_{A}'(z)_{n}) & 0 \\ 0 & \diag(y_{A}'(z)_{1},\dots,y_{A}'(z)_{n})\end{array}\right),\]
where $y_{A}'(z)_{i}=z$ if $v(q_{i})>0$ and $y_{A}'(z)_{i}=\overline{z}$ if $v(q_{i})<0$. We let $y=y'^{\ad}:\mathbb{S}\to H_{v}^{\ad}$. Using the explicit computation of $\gamma^{-1}Q'$ in each case, the group $\operatorname{GU}(\gamma^{-1}Q')$ is isomorphic to a similitude unitary group $\operatorname{GU}(p,q)$ of a certain signature $(p,q)$ (furthermore, if $m=2$, in our case where $D_{0,v}$ is split, the signature is always $(n,n)$, so the group $H_{v}$ is in fact quasi-split). It's then standard that $y'$, and hence $y$, satisfies Deligne's axioms (see for instance the Appendix of \cite{milneshih}). 
\end{rem}

\subsection{Groups of Type $D_{l}$ ($l\geq 4$ odd)}\label{subsec:groupsDl}
Suppose that $\mathcal{D}=D_{l}$ with $l\geq 5$ odd. Then $H=\mathbf{Spin}(B,\tau)$, where $B$ is a central simple algebra over $F$ of degree $2l$ and $\tau$ is an orthogonal involution (\cite{involutions}, Theorem 26.15). The adjoint group is $H^{\ad}=\mathbf{PGO}^{+}(B,\tau)$. In order to avoid introducing spin groups, we will work in this section with $H^{\ad}$. Since the map $\Aut(H)\to\Aut(H^{\ad})$ is an isomorphism, an opposition involution on $H^{\ad}$ will uniquely lift to an opposition involution on $H$; moreover, suppose that $S\subset H$ is a maximal torus and the involution on $H^{\ad}$ preserves $S^{\ad}$ and induces $c_{S_{v}^{\ad}}$ on $S_{v}^{\ad}$ for every $v\in I_{nc}$. Then the lifted involution on $H$ preserves $S$ and also obviously induces $c_{S_{v}^{\ad}}$ on $S_{v}^{\ad}$ for every $v\in I_{nc}$. This will allow us to concentrate on $H^{\ad}$ and avoid spin groups. 

Since $F$ is a number field, it can be shown that $B=\End_{D}(\Lambda)$, where $D=F$ or a quaternion division algebra over $F$ (see \cite{scharlau}, 8.2.3), and $\Lambda$ is a right $D$-vector space of finite dimension $n$. Let $m=\deg_{F}D$. Moreover, the involution $\tau:B\to B$ must be attached to a non-degenerate $F$-bilinear form $q:\Lambda\times\Lambda\to D$. In the case $D=F$ (where $\dim_{F}\Lambda=2l$), $q$ is a symmetric bilinear form. In the case that $D$ is a quaternion division algebra (where $\dim_{D}\Lambda=l$), $q$ is a skew-hermitian form with respect to the canonical involution $\sigma:D\to D$. We will only treat the case where $D$ is a quaternion division algebra. Let $I_{s}\subset I=\Hom(F,\R)$ be the set of $v:F\to\R$ such that $D_{v}=D\otimes_{F,v}\R$ is split, and let $I_{ns}$ be its complement in $I$. For $v\in I_{s}$, the skew-hermitian form $q_{v}$ on $\Lambda_{v}$ defines a non-degenerate symmetric bilinear form $b_{v}$ over a real vector space $W_{v}$ of dimension $2n$ (see \cite{scharlau}), and then we have that $I_{c}\subset I_{s}$ is the set of split places where $b_{v}$ is definite. As in the Appendix of \cite{milneshih} (type $D^{\mathbb{H}}$), we will assume that $I_{c}=I_{s}$. We call the pair $(\Lambda,q)$ a skew-hermitian space over $D$. Note that $n=l$ is odd.

Let $\beta=\{v_{1},\dots,v_{n}\}$ be a $D$-basis of $\Lambda$, which is $q$-orthogonal. The group $H^{\ad}=\mathbf{PGO}^{+}(\Lambda,q)$ can also be seen as the adjoint group of $G=\mathbf{SO}(\Lambda,q)$, where
\[ G(A)=\{g\in\Aut_{D\otimes_{F}{A}}(\Lambda\otimes_{F}{A}):q_{A}(g(x),g(y))=q_{A}(x,y)\quad\forall x,y\in\Lambda_{A},\quad\operatorname{Nrd}(g)=1\}\]
for an $F$-algebra $A$. Here $\operatorname{Nrd}$ is the reduced norm in $\End_{D}(\Lambda)$. We let $S'=S'_{\beta}\subset G$ be the subgroup of $G$ defined as follows. For every $i=1,\dots,n$, let $q_{i}=q(v_{i},v_{i})$. This is a pure quaternion in $D$, and so $L_{i}=F(q_{i})$ is a quadratic field extension of $F$. For an $F$-algebra $A$, let
\[ S'(A)=\{g\in G(A)\subset\Aut_{D\otimes_{F}A}(\Lambda\otimes_{F}A):g(v_{i}\otimes1)=(v_{i}\otimes 1)\lambda_{i},\text{ for some }\lambda_{i}\in(L_{i}\otimes_{F}A)^{\times}\}.\]
Then $S'\subset G$ is a maximal torus of $G$, and it defines maximal tori $S\subset H$ and $S^{\ad}=S'^{\ad}\subset H^{\ad}$.

We will construct involutions on $H$ modeled after our constructions for the case of type $A_{l}$. For this we need to make an analogous extra assumption.

\begin{dfn}\label{dfn:strskewhermitian} We say that the skew-hermitian space $(\Lambda,q)$ over $D$ is {\em strongly skew-hermitian} if there exists a $q$-orthogonal $D$-basis $\beta=\{v_{1},\dots,v_{n}\}$ of $\Lambda$ and an $F$-automorphism $\alpha:D\to D$ such that $q(v_{i},v_{j})=-\alpha(q(v_{j},v_{i}))$ and $\alpha^{2}=1$.
\end{dfn}

\begin{rem} Any automorphism $\alpha:D\to D$ as above must be necessarily inner, of the form $\alpha(d)=rdr^{-1}$ for some $r\in D^{\times}$ such that $r^{2}\in F^{\times}$. This implies that $r\sigma(r)^{-1}\in F^{\times}$ as well (because $F$ is the set of elements of $D$ fixed by $\sigma$). Moreover, since $q(v_{i},v_{i})\in D^{\times}$, $r$ must be a pure quaternion in $D$.
\end{rem}

Suppose that $(\Lambda,q)$ is strongly hermitian, and let $\beta$ and $\alpha=\inte(r)$ be as in the definition. We then have $\alpha\sigma=\sigma\alpha$. Let $I:\Lambda\to\Lambda$ be the $\alpha$-semilinear automorphism obtained by applying $\alpha$ to the coefficients of elements of $\Lambda$ with respect to the basis $\beta$. Then $q(I(x),I(y))=-\alpha(q(x,y))$. Let $\theta:G\to G$ be defied by $\theta_{A}(g)=I_{A}gI_{A}$ for an $F$-algebra $A$ and a $D\otimes_{F}A$-linear automorphism $g$ of $\Lambda\otimes_{F}A$. 

Let $L=F(r)$, where $r\in D$ is as above. This is again a quadratic extension of $F$ (and a maximal subfield of $D$). Let $S'$ and $S^{\ad}$ be the maximal tori of $G$ and $H^{\ad}$ defined above using the basis $\beta$.

\begin{prop} With the above hypotheses, the following statements are true.
\begin{enumerate}[(a)]
\item The map $\theta:G\to G$ is an opposition involution (and hence so is $\theta^{\ad}$).
\item We have $\theta(S')=S'$ and for every $v\in I_{nc}$, $\theta_{v}:G_{v}\to G_{v}$ induces $c_{S'_{v}}$ on $S'_{v}$.
\end{enumerate}
In particular, $\theta^{\ad}_{v}:H^{\ad}_{v}\to H^{\ad}_{v}$ induces $c_{S^{\ad}_{v}}$ on $S^{\ad}_{v}$ for $v\in I_{nc}$.
\begin{proof}
For part (a), it suffices to see that $\theta_{E}:G_{E}\to G_{E}$ is an opposition involution, for a convenient extension $E/F$. Using the basis $\beta$ and the isomorphism $\phi:D\otimes_{F}L\to M_{2}(L)$ as constructed above, we can identify $G_{L}$ as follows. Implicit in the construction of $\phi$ is the choice of a pure quaternion $s\in D$ with $rs=-rs$, and we let $t=v(s^{2})\in\R$. Let $q_{i}=q(v_{i},v_{i})$. Since $\sigma(q_{i})=-q_{i}$ and $rq_{i}r^{-1}=-q_{i}$, we have
\[ \phi(r)=\left(\begin{array}{cc}r & 0 \\ 0 & -r\end{array}\right)\in\GL_{2}(L) \]
and
\[ \phi(q_{i})=\left(\begin{array}{cc}0 & b_{i}\\c_{i}& 0\end{array}\right) \]
for some $b_{i},c_{i}\in L$. The image in $M_{2}(L)$ under $\phi$ of $L_{i}\otimes _{F}L\subset D\otimes_{F}L$ consist of the matrices in $M_{2}(L)$ of the form 
\[ \left(\begin{array}{cc} x & yb_{i} \\ yc_{i}& x\end{array}\right)\]
for some $x,y\in L$. Thus, the induced isomorphism $\phi:M_{n}(D\otimes_{F}L)\to M_{2n}(L)$ sends the subalgebra of diagonal matrices $L_{1}\otimes_{F}L\times\dots\times L_{n}\otimes_{F}L$ to the set of matrices in $M_{2n}(L)$ of the form 
\begin{equation}\label{descrS} X=\left(\begin{array}{cc}\diag(x_{1},\dots,x_{n}) & \diag(y_{1}b_{1},\dots,y_{n}b_{n}) \\ \diag(y_{1}c_{1},\dots,y_{n}c_{n}) & \diag(x_{1},\dots,x_{n})\end{array}\right) \end{equation}
with $x_{i},y_{i}\in L$. Let
\[ \widetilde Q=\left(\begin{array}{cc}0 & \diag(b_{1},\dots,b_{n}) \\ \diag(c_{1},\dots,c_{n}) & 0\end{array}\right)\in\GL_{2n}(L).\]
Then, for any $L$-algebra $R$, writing a matrix $X\in\GL_{2n}(R)$ as $X=\left(\begin{array}{cc}A & B \\ C & D\end{array}\right)$, there is an isomorphism
\begin{equation}\label{identGR} G(R)\cong\left\{X\in\GL_{2n}(R):\left(\begin{array}{cc}{}^{t}D & -{}^{t}B \\ -{}^{t}C & {}^{t}A\end{array}\right)\widetilde Q\left(\begin{array}{cc}A & B\\C & D\end{array}\right)=\widetilde Q,\quad\det(X)=1\right\}\end{equation}
that takes the subgroup $S'$ to the subgroup of matrices of the form (\ref{descrS}) in the right hand side. Note that the equation is equivalent to ${}^{t}X\widetilde Q'X=\widetilde Q'$, where 
\[ \widetilde Q'=\left(\begin{array}{cc}\diag(c_{1},\dots,c_{n}) & 0 \\ 0 & -\diag(b_{1},\dots,b_{n})\end{array}\right)\]
(the matrix $\widetilde Q'$ is the matrix of the associated bilinear form; see \cite{scharlau}, 10.3). Moreover, if
\[ \gamma=\left(\begin{array}{cc}r I_{n} & 0 \\ 0&-rI_{n}\end{array}\right),\]
then $\theta_{R}(X)=\gamma X\gamma^{-1}$ for $X\in G(R)$; in block matrix terms,
\[ \theta_{R}\left(\begin{array}{cc}A & B\\ C & D\end{array}\right)=\left(\begin{array}{cc}A & -B\\ -C & D\end{array}\right).\]
It's clear then that $\theta$ preserves $S'$. 

Let $E/L$ be a field extension such that there exist elements $e_{i},f_{i}\in E$ with $e_{i}^{2}=c_{i}$ and $f_{i}^{2}=b_{i}$ (for example, take $E=\C$ with a fixed embedding of $L$). For elements $a_{1},\dots,a_{n}$, let $\operatorname{adiag}(a_{1},\dots,a_{n})$ be the anti-diagonal matrix whose $(i,n+1-i)$-th entry is $a_{i}$, and let $J_{n}=\operatorname{adiag}(1,\dots,1)$. Let
\[ \delta=\left(\begin{array}{cc}\operatorname{adiag}(e_{n},\dots,e_{1}) & \operatorname{adiag}(-f_{n},\dots,-f_{1})\\ \diag(\frac{e_{1}}{2},\dots,\frac{e_{n}}{2}) & \diag(\frac{f_{1}}{2},\dots,\frac{f_{n}}{2})\end{array}\right)\in\GL_{2n}(E).\]
Then the map $X\mapsto\delta X\delta^{-1}$ sends $G_{E}$ (viewed inside $\GL_{2n,E}$ via (\ref{identGR})) to the special orthogonal group $\SO_{2n}$ of the matrix $J_{2n}$ over $E$. The maximal torus $S'_{E}$ maps to the subgroup of diagonal matrices in $\SO_{2n}$, and $\theta$ becomes conjugation by the matrix 
\[ \delta\gamma\delta^{-1}=\left(\begin{array}{cc} 0 & 2rJ_{n}\\ \frac{r}{2}J_{n} & 0\end{array}\right)\]
inside $\GL_{2n}$. We identify in the usual way $X^{*}(S')\cong\Z^{n}$. As a Borel subgroup of $G_{E}$ we take the subgroup $B$ of upper triangular matrices belonging to $G_{E}$. The map $\theta$ sends $B$ to the subgroup $B^{-}$ of lower triangular matrices. Let $J_{2n}'$ be the matrix obtained from $J_{2n}$ by swapping the rows $n$ and $n+1$. Then it's easy to see that $J_{2n}'\in G(E)$ and sends $B^{-}$ to $B$. It follows that $\Psi_{0}(\theta)(\chi)=\chi\circ\inte(J_{2n}')\circ\theta$ for $\chi\in X^{*}(S')$. If $\chi$ is parametrized by $(a_{1},\dots,a_{n})$, then $\Psi_{0}(\theta)(\chi)$ is parametrized by $(a_{1},\dots,a_{n-1},-a_{n})=(a_{1},\dots,a_{n})^{\star}$ (see \cite{bourbaki}, Plate IV). Thus, $\theta:G\to G$ is an opposition involution.

Let $v:F\hookrightarrow\R$ and let $\tau:L\hookrightarrow\C$ be an extension of $v$ to $L$. If $\tau=\overline\tau$, then $\tau(r)\in\R$. Thus, $\tau(r)^{2}\in\R_{>0}$, and this implies that $D_{v}$ is split, so $v\in I_{s}=I_{c}$. In part (b), we only care for $v\in I_{nc}$, so suppose from now on that $\tau\neq\overline\tau$, so that $\tau(r)\in i\R_{>0}$. By the same reasoning we have that $t=v(s^{2})<0$. We use $\tau$ to identify $G_{\C}\cong\SO_{2n}$ as above. We first work out the induced complex conjugation on $G(\C)\cong\SO_{2n}(\C)$.  Using the isomorphisms $D\otimes_{F,v}\C\simeq(D\otimes_{F}L)\otimes_{L,\tau}\C\cong M_{2}(\C)$ (the last one coming from $\phi$), it's easy to see that complex conjugation on $D\otimes_{F,v}\C$ corresponds to taking a matrix $X\in M_{2}(\C)$ to
\[ \left(\begin{array}{cc}t & 0 \\ 0 & 1\end{array}\right)\left(\begin{array}{cc}\overline{X_{22}} & \overline{X_{21}} \\ \overline{X_{12}} & \overline{X_{11}}\end{array}\right)\left(\begin{array}{cc}t^{-1} & 0 \\ 0 & 1\end{array}\right),\]
where $t=v(s^{2})$ as above. Note that $q_{i}\in D\subset D\otimes_{F,v}\C$, so this implies that $t\overline{\tau(c_{i})}=\tau(b_{i})$ and thus
\begin{equation}\label{eqn:tci}t\frac{e_{i}}{\overline{f_{i}}}=-\frac{f_{i}}{\overline{e_{i}}}.\end{equation} It follows that the induced complex conjugation on $G(\C)$, viewed inside $\GL_{2n}(\C)$ as in (\ref{identGR}), is given by
\[ X=\left(\begin{array}{cc}A&B\\C&D\end{array}\right)\mapsto c'(X)=\left(\begin{array}{cc}\overline{D} & t\overline{C} \\ t^{-1}\overline{B} & \overline{A}\end{array}\right).\]
Finally, we apply conjugation by $\delta$ to identify $G_{\C}$ with $\SO_{2n}$. We only need to consider the action of $c$ on diagonal matrices. Let $X=\diag(x_{1},\dots,x_{n},x_{n}^{-1},\dots,x_{1}^{-1})\in\SO_{2n}(\C)$. Then $c(X)=\delta c'(\delta)^{-1}c'(X)c'(\delta)\delta^{-1}$, and a long but easy direct calculation using (\ref{eqn:tci}) shows that
\[ \delta c'(\delta)^{-1}=\left(\begin{array}{cc}2\operatorname{adiag}(\frac{e_{n}}{\overline{f_{n}}},\dots,\frac{e_{1}}{\overline{f_{1}}}) & 0 \\ 0 & \frac{1}{2}\operatorname{adiag}(\frac{f_{1}}{\overline{e_{1}}},\dots,\frac{f_{n}}{\overline{e_{n}}})\end{array}\right),\]
and thus
\[ c(\diag(x_{1},\dots,x_{n},x_{n}^{-1},\dots,x_{1}^{-1})=\diag(\overline{x_{1}}^{-1},\dots,\overline{x_{n}}^{-1},\overline{x_{n}},\dots,\overline{x_{1}}).\]
This implies that if $\chi\in X^{*}(S')$ is parametrized by $(a_{1},\dots,a_{n})\in\Z^{n}$ then $\chi^{c}$ is parametrized by $(-a_{1},\dots,-a_{n})$. This is also easily seen to be the parameter of $\chi\circ\theta$, which shows that $\theta_{v}$ induces $c_{S_{v}'}$ on $S_{v}'$.

\end{proof}
\end{prop}

\begin{rem}\label{rem:mapyD} Keep the assumptions and notation as above. For each $v\in I_{nc}$, we will construct a map $y:\mathbb{S}\to H_{v}^{\ad}$ satisfying Deligne's axioms (\cite{delignecorvallis}, 1.2.1) and factoring through $S_{v}^{\ad}$. Recall that $t=v(s^{2})$ and let $u=v(r^{2})$. Since $v\in I_{nc}$, by our assumptions $D_{v}$ is not split. This implies that $u<0$ and $t<0$. Let $\psi:D_{v}\to\mathbb{H}$ be the isomorphism of $\R$-algebras sending $r\otimes1$ to $\sqrt{-u}e_{2}$ and $s\otimes 1$ to $\sqrt{-t}e_{3}$. Here $e_{1}$, $e_{2}$, $e_{3}$ and $e_{4}$ are the following elements of $\mathbb{H}$:
\[ e_{1}=I_{2},\quad e_{2}=\left(\begin{array}{cc}i & 0 \\ 0 & -i\end{array}\right),\]
\[ e_{3}=\left(\begin{array}{cc}0 & 1 \\ -1 & 0\end{array}\right),\quad e_{4}=e_{2}e_{3}.\]
As above, we can write $\psi(q_{i})=\left(\begin{array}{cc} 0 & y_{i}\\ -\overline{y_{i}} & 0\end{array}\right)$ with $y_{i}\in\C^{\times}$. Let
\[ T=\left(\begin{array}{cc} 0 & \diag(y_{1},\dots,y_{n}) \\ -\diag(\overline{y_{1}},\dots,\overline{y_{n}})\end{array}\right).\]
We then have, for an $\R$-algebra $R$,
\begin{equation}\label{identGvRH} G_{v}(R)\cong\left\{X=\left(\begin{array}{cc}A & B \\ -\overline{B} & \overline{A}\end{array}\right)\in\GL_{2n}(\C\otimes_{\R}R):X^{*}TX=T,\quad \det(X)=1\right\}.\end{equation}
The maximal torus $S'$ corresponds to the subgroup of matrices on the right hand side where $A=\diag(a_{1},\dots,a_{n})$ and $B=\diag(b_{1}y_{1},\dots,b_{n}y_{n})$ with $a_{i},b_{i}\in R$. We can actually see $H_{v}$ as the adjoint group of $G'_{v}$, where
\[G_{v}'(R)\cong\left\{X=\left(\begin{array}{cc}A & B \\ -\overline{B} & \overline{A}\end{array}\right)\in\GL_{2n}(\C\otimes_{\R}R):X^{*}TX=\nu(X)T,\quad \det(X)=\nu(X)^{n}\right\}.\]
We define $y':\mathbb{S}\to G_{v}'$ by the formula
\[ y'_{R}(z)=\left(\begin{array}{cc}\operatorname{Re}(z)I_{n} & \diag(\frac{\operatorname{Im}(z)}{|y_{1}|}y_{1},\dots,\frac{\operatorname{Im}(z)}{|y_{n}|}y_{n}) \\ \diag(-\frac{\operatorname{Im}(z)}{|y_{1}|}\overline{y_{i}},\dots,-\frac{\operatorname{Im}(z)}{|y_{n}|}\overline{y_{n}}) & \operatorname{Re}(z)I_{n}\end{array}\right)\]
for $z\in\mathbb{S}(R)$. Conjugating by a suitable matrix $U\in\GL_{2n}(\C)$, we can write $G'_{v}\cong\operatorname{GO}^{*}(2n)$ and $y$ becomes the map in the Appendix of \cite{milneshih}, so it satisfies Deligne's axioms, and hence also does $y=y'^{\ad}$.
\end{rem}

\section{Involutions on certain Shimura varieties}\label{sec:invoshimura}

In this section we combine all our previous results to prove the existence of descent data on certain Shimura varieties $\Sh(G,X)$. As we said before, we only consider the case where the simple groups $H_{i}$ are of type $A$ or $D^{\mathbb{H}}$. In the previous section, we constructed opposition involutions on some of these groups, preserving a certain maximal torus $S_{i}$ and inducing complex conjugation on its characters. Furthermore, we constructed maps $y_{i,v}:\mathbb{S}\to H_{i,v}^{\ad}$ for every $v\in I_{i,nc}$ satisfying Deligne's axioms (\cite{delignecorvallis}, 1.2.1), factoring through $S_{i,v}^{\ad}$. We now show that we can always find an element $x\in X$ such that $x_{i,v}$ factors through $S_{i,v}^{\ad}$ for every $i$ and $v\in I_{i,nc}$. The existence of descent data will follow by combining this with Proposition \ref{mainprop}.

Let $H$ be an almost simple, simply connected group over $\R$ (to play the role of one of the non-compact $H_{i,v}$). Suppose that there exist morphisms $y:\mathbb{S}\to H^{\ad}$ satisfying Deligne's axioms (\cite[1.2.1]{delignecorvallis}); in particular, $H$ is absolutely almost simple. Let $D$ be the Dynkin diagram of $H_{\C}$ associated with a choice of maximal torus and Borel. To each $H^{\ad}(\R)$-conjugacy class $Y$ of morphisms $y$ as above, we can attach a special node $s_{Y}\in D$, and $s_{Y}=s_{Y'}$ if and only if $Y=Y'$.

\begin{lemma} Under the above conditions, there exist at most two $H(\R)$-conjugacy classes $Y$ of morphisms satisfying 1.2.1 of \cite{delignecorvallis}. Moreover, given such a conjugacy class $Y$, any morphism satisfying these axioms must belong to either $Y$ or $Y^{-1}$.
\begin{proof}
Suppose first that $D$ is not of type $A_{l}$. This case is easy because there are not too many special nodes. Indeed, assume first that $H(\R)$ is connected, and fix $Y$ one of the conjugacy classes. Then $s_{Y^{-1}}=s_{Y}^{\star}\neq s_{Y}$ (\cite[1.2.8]{delignecorvallis}), and hence $Y^{-1}$ and $Y$ are two distinct conjugacy classes. Suppose that $Z$ is a third conjugacy class, that is, $s_{Z}$ is neither equal to $s_{Y}$ nor to $s_{Y}^{\star}$. Again by \cite[1.2.8]{delignecorvallis}, $s_{Z}\neq s_{Z}^{\star}$, and thus we have four distinct special nodes $s_{Y}$, $s_{Y}^{\star}$, $s_{Z}$ and $s_{Z}^{\star}$. There is no connected Dynkin diagram with four special nodes which is not of type $A_{l}$, and thus this is a contradiction. If $H(\R)$ is not connected, then $s_{Y}=s_{Y}^{\star}$ by {\em op. cit.}. If $Z$ is another conjugacy class, then again by {\em op. cit.} we must have $s_{Z}=s_{Z}^{\star}$. But for any connected Dynkin diagram, there is at most one special node which is fixed under the opposition involution, and thus $Z=Y$. 

Suppose now that $H$ is of type $A_{l}$ with $l\geq 2$, so $H={SU}(p,q)$ for some non-zero pair of integers $p,q$ such that $p+q=l+1$. The isomorphism $\C\otimes_{\R}\C\simeq\C\times\C$ given by $z\otimes a\mapsto(za,\overline za)$ induces by projection on the first coordinate an isomorphism $H_{\C}\simeq\operatorname{SL}_{l+1,\C}$; fix the usual Borel pair here to define the Dynkin diagram. Define a morphism $y_{0}:\mathbb{S}\to H^{\ad}=\operatorname{PGU}(p,q)$ with $y_{0}(z)$ being the class of the matrix \[ \left(\begin{array}{cc} z I_{p} & 0 \\ 0 & \overline z I_{q}\end{array}\right).\] Then $y_{0}$ satisfies axioms 1.2.1 of \cite{delignecorvallis}, and the special node $s_{0}$ attached to its $H^{\ad}(\R)$-conjugacy class $Y_{0}$ is $\alpha_{p}$. From the conjugate map $\overline y_{0}=y_{0}^{-1}$ we get the special node $\alpha_{q}$ associated with $Y_{0}^{-1}$. If $Y$ is another conjugacy class, say with special node $\alpha_{t}$, then there would be an isomorphism $\operatorname{PGU}(p,q)\cong\operatorname{PGU}(t,l+1-t)$ sending $Y_{0}$ or $Y_{0}^{-1}$ to $Y$. In particular, $t=p$ or $t=q$, and we conclude that there are at most two possible conjugacy classes of morphisms satisfying 1.2.1 of \cite{delignecorvallis} for the fixed form $\operatorname{PGU}(p,q)$ of $\operatorname{PGL}_{l+1,\C}$ (and there are exactly two in all cases except when $p=q$, when there is only one).
\end{proof}
\end{lemma}

Going back to our general Shimura datum $(G,X)$, for each $i$, let $S_{i}\subset H_{i}$ be a maximal torus, $\widetilde T_{i}=\Res_{F_{i}/\Q}S_{i}\subset\widetilde G_{i}$, $T_{i}\subset G_{i}$ its image in $G_{i}$, $T'\subset G^{\der}$ the image of their product, and $T=Z^{0}T'$. Note that $T_{\R}^{\ad}=T_{\R}^{'\ad}=\prod_{i,v}S_{i,v}^{\ad}$, where $S_{i,v}\subset H_{i,v}$ and $S_{i,v}^{\ad}$ is its image in $H_{i,v}^{\ad}$. 

\begin{lemma}\label{lemmatorus} Suppose that $T\subset G$ is a the maximal torus defined above. Suppose that for each $v\in I_{i,nc}$, there exists a morphism $y_{i,v}:\mathbb{S}\to H_{i,v}^{\ad}$ satisfying axioms 1.2.1 of \cite{delignecorvallis} and factoring through $S_{i,v}^{\ad}$. Then there exists an element $x\in X$ such that $x^{\ad}$ factors through $T_{\R}^{\ad}$.
\begin{proof} Let $z\in X$ be an arbitrary element. The previous lemma implies that $z_{i,v}$ is $H_{i,v}^{\ad}(\R)$-conjugate to a map $y_{i,v}:\mathbb{S}\to S_{i,v}^{\ad}$. Thus, we can write $z_{i,v}=u_{i,v}.y_{i,v}$ for $u_{i,v}\in H_{i,v}^{\ad}(\R)$. We claim that, after possibly changing the $y_{i,v}$, we can arrange for $u_{i,v}$ to be in $H_{i,v}^{\ad}(\R)^{+}$. Indeed, if $u_{i,v}$ is not in that connected component, then in particular $H_{i,v}^{\ad}(\R)$ is not connected, and thus there is only one conjugacy class in question, with two connected components, one containing $z_{i,v}$ and the other one containing $y_{i,v}$. Thus, we only need to replace $y_{i,v}$ with $y_{i,v}^{-1}$, which also factors through $S_{i,v}^{\ad}$. For $v\in I_{i,c}$, let $u_{i,v}=1$. It follows that $u=(u_{i,v})\in G^{\ad}(\R)^{+}$, and thus there exists $g\in G(\R)$ lifting $u$. Let $x=g^{-1}.z\in X$, so that $x^{\ad}=(y_{i,v})$, which factors through $T_{\R}^{\ad}$ as desired.
\end{proof}
\end{lemma}

\begin{dfn} The Shimura datum $(G,X)$ is said to be {\em strongly of type} $(AD^{\mathbb{H}})$ if each of the groups $H_{i}$ is either of type $A_{l}$ with $l\geq 2$ and attached to a strongly hermitian space (as in Definition \ref{dfn:strhermitian}), or of type $D_{l}$ with $l\geq 5$ odd and attached to a strongly skew-hermitian space (as in Definition \ref{dfn:strskewhermitian}).
\end{dfn}

For example, a Shimura variety defined by a similitude unitary group attached to a hermitian space over a CM field is strongly of type $(AD^{\mathbb{H}})$. Note however that the definition only restricts the semisimple part of $G$.

\begin{thm}\label{thmforstrong} Suppose that $(G,X)$ is strongly of type $(AD^{\mathbb{H}})$. Then there exists an involution $\theta:G\to G$ such that $\theta(X)=\overline{X}$, and hence there exist a model of $\Sh(G,X)$ over $E^{+}$ as in Theorem \ref{thmdescent}.
\begin{proof} In Subsections \ref{subsec:groupsAl} and \ref{subsec:groupsDl}, we constructed for every $i$, an opposition involution $\theta_{i}:H_{i}\to H_{i}$ and a maximal torus $S_{i}\subset H_{i}$ such that $\theta_{i}(S_{i})=S_{i}$ and $\theta_{i,v}^{\ad}$ induces $c_{S_{i,v}^{\ad}}$ for every $v\in I_{i,nc}$. Moreover, by Remarks \ref{rem:mapyA} and \ref{rem:mapyD}, for every $i$ and $v\in I_{i,nc}$, there is a map $y_{i,v}:\mathbb{S}\to H_{i,v}^{\ad}$ satisfying Deligne's axioms (\cite{delignecorvallis}, 1.2.1) and factoring through $S_{i,v}^{\ad}$. The result then follows by combining Proposition \ref{mainprop} and Lemma \ref{lemmatorus}.

\end{proof}
\end{thm}

\begin{rem} The conclusion of the previous theorem holds in other cases as well. For instance, if $G$ is adjoint and there exists an opposition involution $\theta:G\to G$ (which is always the case if $G$ is also quasi-split, for example), then by the adjointness of $G$, we conclude that $\theta(X)=\overline{X}$. Our method can also work to include factors of other types, for instance of type $E_{6}$, as long as one can construct an opposition involution inducing complex conjugation on the characters of a maximal torus (at non-compact places) and morphisms satisfying Deligne's axioms and factoring through these tori (as in Remarks \ref{rem:mapyA} and \ref{rem:mapyD}).
\end{rem}

\bibliography{realmodels}{}
\bibliographystyle{alpha}

\vspace{3mm}
\noindent{Don Blasius\\
Mathematics Department, University of California, Los Angeles, CA 90024, USA\\
Email: {\tt blasius@math.ucla.edu}\\
Webpage: {\tt http://www.math.ucla.edu/\~{}blasius/}}

\vspace{3mm}
\noindent{Lucio Guerberoff\\
Mathematics Department, University College London, 25 Gordon Street, London WC1H 0AY, UK\\
Email: {\tt l.guerberoff@ucl.ac.uk}\\
Webpage: {\tt http://www.ucl.ac.uk/\~{}ucahlgu/}}

\end{document}